\documentclass[12pt]{amsart}
\usepackage{amsmath,amscd,amssymb}
\usepackage[all]{xy}
\theoremstyle{plain}
\newtheorem{theorem}{Theorem}[section]
\newtheorem{proposition}[theorem]{Proposition}
\newtheorem{lemma}[theorem]{Lemma}
\newtheorem{corollary}[theorem]{Corollary}
\theoremstyle{definition}
\newtheorem{definition}[theorem]{Definition}
\newtheorem{example}[theorem]{Example}
\theoremstyle{remark}
\newtheorem{remark}[theorem]{Remark}
\numberwithin{equation}{section}
\newcommand{\benu}{\begin{enumerate}\renewcommand{\labelenumi}{{\rm (\roman{enumi})}}\renewcommand{\itemsep}{0pt}}
\newcommand{\eenu}{\end{enumerate}}
\newcommand{\N}{\mathbb{N}}

\newcommand{\R}{\mathbb{R}}
\newcommand{\C}{\mathbb{C}}
\newcommand{\T}{\mathbb{T}}
\newcommand{\e}{\varepsilon}
\newcommand{\cK}{{\mathcal K}}
\newcommand{\cL}{{\mathcal L}}

\newcommand{\cO}{{\mathcal O}}
\newcommand{\cT}{{\mathcal T}}

\newcommand{\ip}[2]{\langle{#1},{#2}\rangle}
\newcommand{\lip}[2]{{}_X\langle{#1},{#2}\rangle}
\newcommand{\bip}[2]{\big\langle{#1},{#2}\big\rangle}
\newcommand{\ti}[1]{\widetilde{{#1}}}
\newcommand{\Ca}{$C^*$-al\-ge\-bra }
\newcommand{\Cas}{$C^*$-al\-ge\-bras }
\newcommand{\Cc}{$C^*$-cor\-re\-spon\-dence }
\newcommand{\Ccs}{$C^*$-cor\-re\-spon\-dences }
\newcommand{\Csa}{$C^*$-sub\-al\-ge\-bra }
\newcommand{\Csas}{$C^*$-sub\-al\-ge\-bras }
\DeclareMathOperator{\id}{id}
\DeclareMathOperator{\Aut}{Aut}

\DeclareMathOperator{\cspa}{\overline{span}}
\DeclareMathOperator{\diag}{diag}

\begin{document}
\title[Ideal structure of $C^*$-algebras 
associated with $C^*$-correspondences]
{Ideal structure of \boldmath{$C^*$}-algebras 
associated with \boldmath{$C^*$}-correspondences}
\author{Takeshi Katsura}
\address{
Department of Mathematical Sciences,
University of Tokyo, Komaba, Tokyo, 153-8914, JAPAN}
\curraddr{
Department of Mathematics,
University of Oregon, 
Eugene, Oregon, 97403-1222, U.S.A.}
\email{katsu@ms.u-tokyo.ac.jp}
\thanks{The author was supported in part by a Research Fellowship 
for Young Scientists of the Japan Society for the Promotion of Science.}



\begin{abstract}
We study the ideal structure of $C^*$-al\-ge\-bras 
arising from $C^*$-cor\-re\-spon\-denc\-es. 
We prove that gauge-invariant ideals of 
our $C^*$-al\-ge\-bras are parameterized 
by certain pairs of ideals of original $C^*$-al\-ge\-bras. 
We show that our $C^*$-al\-ge\-bras have a nice property 
which should be possessed 
by generalization of crossed products. 
Applications to crossed products by Hilbert $C^*$-bimodules 
and relative Cuntz-Pimsner algebras 
are also discussed. 
\end{abstract}

\maketitle

\setcounter{section}{-1}

\section{Introduction}

For a \Ca $A$, 
a \Cc over $A$ is a (right) Hilbert $A$-module 
with a left action of $A$. 
Since endomorphisms (or families of endomorphisms) 
of $A$ define \Ccs over $A$, 
we can regard \Ccs 
as (multi-valued) generalizations of automorphisms or endomorphisms. 
This point of view has same philosophy 
as the idea 
that topological correspondences defined in \cite{Ka2} 
are generalizations of continuous maps 
(see \cite[Section 1]{Ka2}). 

A crossed product by an automorphism 
is a \Ca which has an original \Ca as a $C^*$-sub\-al\-ge\-bra, 
and reflects many aspects of the automorphism. 
For example, the set of ideals of the crossed product 
which are invariant under the dual action 
of the one-dimensional torus $\T$ 
corresponds bijectively to the set of ideals of the original \Ca 
which are invariant under the automorphism. 
As \Ccs are generalizations of endomorphisms, 
a natural problem is to define 
``crossed products'' by $C^*$-cor\-re\-spon\-dences. 
There are plenty of evidence that 
the construction of the \Ca $\cO_X$ from a \Cc $X$ in \cite{Ka4} 
is the right one. 
One piece of evidence is that this 
generalizes many constructions 
which were or were not considered 
as generalizations of crossed products (see \cite{Ka4}). 
We are going to explain another piece of evidence.
For a \Cc $X$, 
we can naturally define a notion of representations of $X$ 
(Definition \ref{DefRep}). 
Thus one \Ca which is naturally associated with a \Cc $X$ 
is a \Ca $\cT_X$ having a universal property 
with respect to representations of $X$ 
(Definition \ref{DefTX}). 
This \Ca $\cT_X$ is nothing 
but an (augmented) Cuntz-Toeplitz algebra 
defined in \cite{Pi}. 
When a \Cc $X$ is defined by an automorphism, 
the \Ca $\cT_X$ is isomorphic 
to the Toeplitz extension 
of the crossed product by the automorphism 
defined in \cite{PV}. 
This \Ca is too large to reflect the informations of $X$. 
In order to get ``crossed products'', 
we have to go to a quotient of $\cT_X$. 
There are two ways to go. 
One way is to define ``covariance'' of representations 
of a $C^*$-cor\-re\-spon\-denc\-e $X$, 
and define a ``crossed product'' by $X$ so that 
it has the universal property 
with respect to covariant representations of $X$. 
This kind of method has been used in many papers, 
and we define our \Ca $\cO_X$ along this line 
(Definition \ref{DefCov}, Definition \ref{DefOX}). 
The other way is to list up the properties of $\cT_X$ 
which the ``crossed product'' should have, 
and define a ``crossed product'' by $X$ 
to be the smallest quotient of $\cT_X$ 
among the quotients satisfying these properties. 
For such properties, 
the following two seem to be reasonable; 
\benu
\item the original \Ca is embedded into the ``crossed product'', 
\item there exists a ``dual action'' of $\T$ on the ``crossed product''. 
\eenu
In this paper, 
we show that these two methods give the same \Ca $\cO_X$ 
(Proposition \ref{smallest}). 
This indicates that the \Ca $\cO_X$ 
is the right one for a ``crossed product'' by a \Cc $X$. 
We note that Cuntz-Pimsner algebras do not satisfy 
the property (i) above 
when the left action of the \Cc is not injective, 
and that the \Ca $\cO_X$ is isomorphic to the Cuntz-Pimsner algebra 
when the left action of the \Cc is injective. 

The ``dual action'' of $\T$ on the \Ca $\cO_X$ 
is called the gauge action. 
The main purpose of this paper 
is to describe the all ideals of the \Ca $\cO_X$ 
associated with a \Cc $X$ 
which are invariant under the gauge action. 
We define invariance of ideals 
of $A$ with respect to a \Cc $X$ over $A$ 
(Definition \ref{Definv}). 
Unlike the case of crossed products by automorphisms, 
we need extra ideals of $A$ other than invariant ideals 
to describe all gauge-invariant ideals of $\cO_X$. 
Similar facts were observed in many papers 
(\cite{BHRS,DT,Ka1,Ka3} to name a few) 
for \Cas arising from graphs or topological graphs. 
We introduce a notion of $O$-pairs 
which are pairs consisting of invariant ideals 
and extra ideals of $A$, 
and show that gauge-invariant ideals 
are parameterized by $O$-pairs (Theorem \ref{GII}). 

This paper is organized as follows. 
In Sections \ref{HilbMod} and \ref{Corres}, 
we fix notations and gather results 
on Hilbert $C^*$-modules 
and $C^*$-cor\-re\-spon\-denc\-es.
In Section \ref{C*algOfCor}, 
we give the definition of our \Cas $\cO_X$ 
constructed from \Ccs $X$. 
In Sections \ref{SecInv} and \ref{SecPairs}, 
we introduce and study invariance of ideals, 
$T$-pairs and $O$-pairs. 
These are related to representations 
of $C^*$-cor\-re\-spon\-denc\-es. 
In Section \ref{SecC*Ad}, 
we construct a \Cc $X_\omega$ from a $T$-pair $\omega$, 
and in Section \ref{SecC^*Rep} 
we prove that this \Cc $X_\omega$ 
has a certain universal property. 
As a corollary, 
we give an alternative definition of our \Cas $\cO_X$ described above 
(Proposition \ref{smallest}). 
In Section \ref{SecIdeal}, 
we prove the main theorem (Theorem \ref{GII}) 
which says that 
the set of all gauge-invariant ideals of $\cO_X$ 
corresponds bijectively to the set of all $O$-pairs of $X$. 
We also see that a quotient of $\cO_X$ by a gauge-invariant ideal 
falls into the class of our $C^*$-al\-ge\-bras. 
In Section \ref{SecGI&Morita}, 
we see that every gauge-invariant ideals have hereditary and full 
\Csas which are isomorphic to \Cas associated with 
$C^*$-cor\-re\-spon\-denc\-es. 
As a consequence of the study 
of crossed products by Hilbert $C^*$-bimodule in Section \ref{SecBimod}, 
all gauge-invariant ideals themselves are shown to be isomorphic to 
\Cas associated with $C^*$-cor\-re\-spon\-denc\-es. 
In Section \ref{SecRCP}, 
we apply our investigation to 
relative Cuntz-Pimsner algebras defined in \cite{MS}.

\medskip

The author is grateful 
to Yasuyuki Kawahigashi for his constant encouragement. 
This paper was written 
while the author was staying at University of Oregon. 
He would like to thank people there for their warm hospitality. 

\medskip

We denote by $\N=\{0,1,2,\ldots\}$ 
the set of natural numbers, 
and by $\C$ the set of complex numbers. 
We denote by $\T$ the group consisting of complex numbers 
whose absolute values are $1$. 
We use a convention that 
$\gamma(A,B)=\{\gamma(a,b)\in D\mid a\in A,b\in B\}$ 
for a map $\gamma\colon A\times B\to D$ 
such as inner products, multiplications or representations. 
We denote by $\cspa\{\cdots\}$ 
the closure of linear spans of $\{\cdots\}$. 
Thus the Cohen factorization theorem can be stated as follows; 

\medskip

\noindent
{\bfseries Lemma.} 
Let $A$ be a $C^*$-al\-ge\-bra, 
$X$ be a Banach space, 
and $\pi\colon A\to B(X)$ be a homomorphism from $A$ 
to the Banach algebra $B(X)$ of the bounded operators on $X$. 
Then we have $\pi(A)X=\cspa (\pi(A)X)$. 

\medskip

We use this result
just to make notation and arguments short. 
The readers who are not familiar with the theorem 
may use $\cspa(\pi(A)X)$ instead of $\pi(A)X$, 
which are actually the same.

\section{Hilbert $C^*$-modules}\label{HilbMod}

\begin{definition}\label{DefHilb}
Let $A$ be a $C^*$-al\-ge\-bra. 
A (right) {\em Hilbert $A$-module} $X$ is a linear space with
a right action of the \Ca $A$ 
and an $A$-valued inner product $\ip{\cdot}{\cdot}_X$ 
satisfying certain conditions 
such that $X$ is complete with respect to the norm defined by 
$\|\xi\|_X=\|\ip{\xi}{\xi}_X\|^{1/2}$ for $\xi\in X$. 
\end{definition}

For a precise definition of Hilbert $C^*$-modules, 
consult \cite{Lnc}. 
We do not assume that 
a Hilbert $A$-module $X$ is full. 
Thus $\cspa\ip{X}{X}_X$ 
can be a proper ideal of $A$, 
where an ideal of a \Ca always means a closed two-sided ideal 
except in the proof of Lemma \ref{LemPed}. 

\begin{definition}
For a Hilbert $A$-module $X$, 
we denote by $\cL(X)$ 
the \Ca of all adjointable operators on $X$. 
For $\xi,\eta\in X$, 
the operator $\theta_{\xi,\eta}\in\cL(X)$ is defined 
by $\theta_{\xi,\eta}(\zeta)=\xi\ip{\eta}{\zeta}_X$ for $\zeta\in X$. 
We define the ideal $\cK(X)$ of $\cL(X)$ by 
$$\cK(X)=\cspa\{\theta_{\xi,\eta}\in \cL(X)\mid \xi,\eta\in X\}.$$
\end{definition}

Let us fix a \Ca $A$ and 
a Hilbert $A$-module $X$ throughout this section. 

\begin{proposition}\label{XI}
Let $I$ be an ideal of $A$. 
For $\xi\in X$, the following are equivalent:
\benu
\item $\xi\in XI$, 
\item $\ip{\eta}{\xi}_X\in I$ for all $\eta\in X$, 
\item $\ip{\xi}{\xi}_X\in I$, 
\item there exist $\eta\in X$ and 
a positive element $a\in I$ such that $\xi=\eta a$. 
\eenu
\end{proposition}

\begin{proof}
Clearly (iv) $\Rightarrow$ (i) $\Rightarrow$ (ii) $\Rightarrow$ (iii). 
For $\xi\in X$ with $\ip{\xi}{\xi}_X\in I$, 
we can find $\eta\in X$ such that $\xi=\eta a$ for 
$a=(\ip{\xi}{\xi}_X)^{1/3}\in I$ (\cite[Lemma 4.4]{Lnc}). 
This proves (iii) $\Rightarrow$ (iv). 
\end{proof}

\begin{corollary}
For an ideal $I$ of $A$, 
$XI$ is a closed linear subspace of $X$ 
which is invariant by the right action of $A$ 
and by the left action of $\cL(X)$.
\end{corollary}

\begin{proof}
Since the set of $\xi\in X$ satisfying the condition (ii) in 
Proposition \ref{XI} is a closed linear space, 
we see that $XI$ is a closed linear space 
(this also follows from the Cohen factorization theorem). 
The rest of the statement is easy to verify. 
\end{proof}

By this corollary, 
$XI$ is a Hilbert $A$-submodule of $X$. 
We can and will consider $\cK(XI)$ 
as a subalgebra of $\cK(X)$ by
$$\cK(XI)=\cspa\{\theta_{\xi,\eta}\in\cK(X) \mid \xi,\eta\in XI\}
\subset \cK(X),$$
(cf. \cite[Lemma 2.6 (1)]{FMR}).
Note that $XI$ is also considered as a Hilbert $I$-module. 
For an ideal $I$ of $A$, 
we denote by $X_I$ the quotient space $X/XI$. 
Both of the natural quotient maps $A\to A/I$ and $X\to X_I$ 
are denoted by $[\cdot]_I$. 
The space $X_I$ has 
an $A/I$-valued inner product $\ip{\cdot}{\cdot}_{X_I}$ 
and a right action of $A/I$ so that 
$$\bip{[\xi]_I}{[\zeta]_I}_{X_I}=\big[\ip{\xi}{\zeta}_X\big]_I,\quad 
[\xi]_I[a]_I=[\xi a]_I$$
for $\xi,\zeta\in X$ and $a\in A$. 
By Proposition \ref{XI}, 
$\eta\in X_I$ satisfies $\ip{\eta}{\eta}_{X_I}=0$ 
only when $\eta=0$. 
Hence $\|\eta\|_{X_I}=\|\ip{\eta}{\eta}_{X_I}\|^{1/2}$ 
defines a norm on $X_I$. 

\begin{lemma}\label{lift}
For $\eta\in X_I$, 
there exists $\xi\in X$ 
such that $\eta=[\xi]_I$ and $\|\eta\|_{X_I}=\|\xi\|_X$. 
\end{lemma}

\begin{proof}
Clearly $[\cdot]_I$ is a norm-decreasing map. 
Thus it suffices to find $\xi\in X$ 
such that $[\xi]_I=\eta$ and $\|\xi\|_X\leq\|\eta\|_{X_I}$
for $\eta\in X_I$. 
Set $C=\|\eta\|_{X_I}^2=\|\ip{\eta}{\eta}_{X_I}\|$. 
Let $f,g$ be functions on $\R_+=[0,\infty)$ defined by 
$$f(r)=\left\{\begin{array}{ll}
1& (0\leq r\leq C)\\
\sqrt{C/r}& (r>C)
\end{array}\right. ,\quad
g(r)=\min\{r,C\}.$$
Then we have $g(r)=rf(r)^2$ and $g(r)\leq C$ for $r\in\R_+$. 
Take $\xi_0\in X$ with $\eta=[\xi_0]_I$. 
Set $a=f(\ip{\xi_0}{\xi_0}_{X})\in \widetilde{A}$ 
and $\xi=\xi_0 a\in X$ where $\widetilde{A}$ 
is the unitization of $A$. 
We have 
$\ip{\xi}{\xi}_X=a^*\ip{\xi_0}{\xi_0}_{X}a=
g(\ip{\xi_0}{\xi_0}_{X})$. 
Hence we get $\|\xi\|_X\leq C^{1/2}=\|\eta\|_{X_I}$. 
Since $f$ is 1 on $[0,C]$, 
we have 
$$[a]_I=f([\ip{\xi_0}{\xi_0}_{X}]_I)=f(\ip{\eta}{\eta}_{X_I})=1.$$ 
Therefore 
we see that $[\xi]_I=[\xi_0]_I[a]_I=\eta$. 
We are done. 
\end{proof}

By this lemma, the norm $\|\cdot\|_{X_I}$ of $X_I$ 
coincides with the quotient norm of $[\cdot]_I\colon X\to X_I$ 
(cf. \cite[Lemma 2.1]{FMR}). 
Hence $X_I$ is complete, and so it is a Hilbert $A/I$-module. 

Since $XI$ is closed under the action of $\cL(X)$, 
we can define a map $\cL(X)\to \cL(X_I)$, 
which is also denoted by $[\cdot]_I$, so that 
$[S]_I[\xi]_I=[S\xi]_I$ for $S\in \cL(X)$ and $\xi\in X$. 
By definition, $S\in\cL(X)$ satisfies $[S]_I=0$ 
if and only if $S\xi\in XI$ for all $\xi\in X$, 
which is equivalent to the condition that 
$\ip{\eta}{S\xi}\in I$ for all $\xi,\eta\in X$ 
by Proposition \ref{XI}. 

\begin{lemma}[{cf. \cite[Lemma 2.6]{FMR}}]\label{ckXI}
For $\xi,\eta\in X$, 
we have $[\theta_{\xi,\eta}]_I=\theta_{[\xi]_I,[\eta]_I}$. 
The restriction of the map $[\cdot]_I\colon \cL(X)\to \cL(X_I)$ 
to $\cK(X)$ is a surjection onto $\cK(X_I)$ 
whose kernel is $\cK(XI)$. 
\end{lemma}

\begin{proof}
The first assertion is easily verified by the definition. 
This implies that the restriction of the map $[\cdot]_I$ 
to $\cK(X)$ is a surjection onto $\cK(X_I)$, 
and that $\cK(XI)$ is in the kernel of $[\cdot]_I$. 
We will show that if $k\in \cK(X)$ satisfies that $[k]_I=0$, 
then $k\in \cK(XI)$. 

There exists an approximate unit 
$\{u_\lambda\}_{\lambda\in\Lambda}$ of $\cK(X)$ 
such that for each $\lambda\in\Lambda$, 
$u_\lambda$ is a finite linear sum of 
elements in the form $\theta_{\xi,\eta}$. 
Take $k\in\cK(X)$ with $[k]_I=0$. 
Since we have $k=\lim ku_\lambda$, 
to prove $k\in \cK(XI)$ 
it suffices to show that $k\theta_{\xi,\eta}\in \cK(XI)$ 
for arbitrary $\xi,\eta\in X$. 
Since $k\xi\in XI$, 
we can find $\xi_0\in X$ and 
a positive element $a_0\in I$
such that $k\xi=\xi_0a_0$ by Proposition \ref{XI}. 
Then we have 
$$k\theta_{\xi,\eta}=\theta_{k\xi,\eta}=\theta_{\xi_0a_0,\eta}
=\theta_{\xi_0\sqrt{a_0},\eta \sqrt{a_0}}\in \cK(XI).$$
We are done. 
\end{proof}

Note that it often happens that $[S]_I\in\cK(X_I)$ 
even if $S\notin\cK(X)$. 
This observation plays an important role 
in our analysis after Section \ref{SecPairs}.
Note also that though three maps $[\cdot]_I\colon A\to A/I$, 
$[\cdot]_I\colon X\to X_I$ and $[\cdot]_I\colon \cK(X)\to \cK(X_I)$ 
are always surjective, 
the map $[\cdot]_I\colon \cL(X)\to \cL(X_I)$ 
need not be surjective 
(because Tietze's extension theorem fails in general). 

Take two ideals $I$ and $I'$ of $A$ such that $I\subset I'$. 
Then $I'/I$ is an ideal of $A/I$ 
and 
$(A/I)/(I'/I)\ni \big[[a]_{I}\big]_{I'/I}\mapsto 
[a]_{I'}\in A/I'$ gives a well-defined isomorphism. 
By this isomorphism, 
we will identify $(A/I)/(I'/I)$ with $A/I'$. 
Thus 
the quotient map $[\cdot]_{I'}\colon A\to A/I'$ coincides with 
the composition of $[\cdot]_{I}\colon A\to A/I$ 
and $[\cdot]_{I'/I}\colon A/I\to A/I'$. 
Similarly we will identify $(X_{I})_{I'/I}$ 
with $X_{I'}$ so that 
$[\cdot]_{I'}=[\cdot]_{I'/I}\circ [\cdot]_{I}$ 
holds for both $X\to X_{I'}$ and $\cL(X)\to\cL(X_{I'})$. 
It is easy to see the following. 

\begin{lemma}\label{isom0}
We have $(XI')_{I}=X_{I}(I'/I)$ in $X_{I}$. 
\end{lemma}

Now take two ideals $I_1$ and $I_2$ of $A$. 
It is well-known that the ideal $I_1\cap I_2$ 
coincides with $I_1I_2$, 
and that $I_1+I_2$ is an ideal of $A$. 
It is easy to see that 
the natural map $I_1/(I_1\cap I_2)\to (I_1+I_2)/I_2$ 
is an isomorphism. 
A pull-back \Ca $B$ of 
two quotient maps $[\cdot]_{(I_1+I_2)/I_1}\colon A/I_1\to A/(I_1+I_2)$ 
and $[\cdot]_{(I_1+I_2)/I_2}\colon A/I_2\to A/(I_1+I_2)$ 
is defined by 
$$B=\big\{(b_1,b_2)\in A/I_1\oplus A/I_2\ \big|\ 
[b_1]_{(I_1+I_2)/I_1}=[b_2]_{(I_1+I_2)/I_2}\in A/(I_1+I_2)\big\}.$$
It is not difficult to see the following 
(see the proof of Proposition \ref{isom3}). 

\begin{lemma}\label{isom1}
The map 
$$\varPi\colon A/(I_1\cap I_2)\ni b\mapsto 
\big([b]_{I_1/(I_1\cap I_2)},[b]_{I_2/(I_1\cap I_2)}\big)\in B$$
is an isomorphism. 
\end{lemma}

We will show analogous statements 
for Hilbert modules and sets of operators on them. 
Define a linear space $Y$ by 
$$Y=\big\{(\eta_1,\eta_2)\in X_{I_1}\oplus X_{I_2}\ \big|\ 
[\eta_1]_{(I_1+I_2)/I_1}=[\eta_2]_{(I_1+I_2)/I_2}\in X_{I_1+I_2}\big\}.$$
We define a $B$-valued inner product on $Y$ by 
$$\bip{(\eta_1,\eta_2)}{(\eta_1',\eta_2')}_Y
=\big(\ip{\eta_1}{\eta_1'}_{X_{I_1}},\ip{\eta_2}{\eta_2'}_{X_{I_2}}\big)
\in B,$$
for $(\eta_1,\eta_2),(\eta_1',\eta_2')\in Y$. 
Clearly $Y$ is complete with respect to the norm 
defined by the inner product. 
If we define a right action of $B$ on $Y$ by 
$$(\eta_1,\eta_2)(b_1,b_2)=(\eta_1b_1,\eta_2b_2)\in Y$$
for $(\eta_1,\eta_2)\in Y, (b_1,b_2)\in B$, 
then we can easily see that $Y$ is a Hilbert $B$-module. 

\begin{lemma}\label{isom2}
The restriction of the quotient map 
$[\cdot]_{I_2/(I_1\cap I_2)}\colon  X_{I_1\cap I_2}\to X_{I_2}$ 
to $X_{I_1\cap I_2}\big(I_1/(I_1\cap I_2)\big)$ 
is a bijection onto 
$X_{I_2}\big((I_1+I_2)/I_2\big)$. 
\end{lemma}

\begin{proof}
By Lemma \ref{isom0}, 
we have 
$X_{I_1\cap I_2}\big(I_1/(I_1\cap I_2)\big)=(XI_1)_{I_1\cap I_2}$. 
It is easy to see that 
the surjection 
$[\cdot]_{I_2/(I_1\cap I_2)}\colon (XI_1)_{I_1\cap I_2}\to (XI_1)_{I_2}$
is injective. 
It is also easy to see that $(XI_1)_{I_2}=\big(X(I_1+I_2)\big)_{I_2}$. 
We have $\big(X(I_1+I_2)\big)_{I_2}=X_{I_2}\big((I_1+I_2)/I_2\big)$ 
by Lemma \ref{isom0}. 
This completes the proof. 
\end{proof}

\begin{proposition}\label{isom3}
By $\varPi$ in Lemma \ref{isom1}, 
we can consider $X_{I_1\cap I_2}$ as a Hilbert $B$-module. 
Then the map 
$$T\colon X_{I_1\cap I_2}\ni\eta\mapsto 
\big([\eta]_{I_1/(I_1\cap I_2)},[\eta]_{I_2/(I_1\cap I_2)}\big)\in Y$$
is an isomorphism as Hilbert $B$-modules. 
\end{proposition}

\begin{proof}
Clearly $T$ preserves inner products 
and right actions. 
This implies that $T$ is isometric. 
It remains to show that $T$ is surjective. 
Take $(\eta_1,\eta_2)\in Y$. 
Since $[\cdot]_{I_1/(I_1\cap I_2)}\colon X_{I_1\cap I_2}\to X_{I_1}$ is 
surjective, 
we can find $\eta'\in X_{I_1\cap I_2}$ 
with $[\eta']_{I_1/(I_1\cap I_2)}=\eta_1$. 
Since 
$[\eta_2]_{(I_1+I_2)/I_2}=[\eta_1]_{(I_1+I_2)/I_1}
=[\eta']_{(I_1+I_2)/(I_1\cap I_2)}$, 
we have 
$$\eta_2-[\eta']_{I_2/(I_1\cap I_2)}\in 
\ker ([\cdot]_{(I_1+I_2)/I_2})=X_{I_2}\big((I_1+I_2)/I_2\big).$$
By Lemma \ref{isom2}, 
we can find $\eta''\in X_{I_1\cap I_2}\big(I_1/(I_1\cap I_2)\big)$ 
with 
$$[\eta'']_{I_2/(I_1\cap I_2)}=\eta_2-[\eta']_{I_2/(I_1\cap I_2)}.$$ 
Set $\eta=\eta'+\eta''\in X_{I_1\cap I_2}$. 
We see that 
\begin{align*}
[\eta]_{I_1/(I_1\cap I_2)}
&=[\eta']_{I_1/(I_1\cap I_2)}+0
=\eta_1,\\
[\eta]_{I_2/(I_1\cap I_2)}
&=[\eta']_{I_2/(I_1\cap I_2)}+[\eta'']_{I_2/(I_1\cap I_2)}
=\eta_2.
\end{align*}
Therefore $T(\eta)=(\eta_1,\eta_2)$. 
Thus $T$ is surjective. 
\end{proof}

\begin{proposition}\label{isom4}
Let us define a \Ca ${\mathcal M}$ by 
$${\mathcal M}=
\big\{(S_1,S_2)\in\cL(X_{I_1})\oplus \cL(X_{I_2})\ \big|\ 
[S_1]_{(I_1+I_2)/I_1}=[S_2]_{(I_1+I_2)/I_2}\in \cL(X_{I_1+I_2})\big\}.$$
Then the map 
$$\varPsi\colon \cL(X_{I_1\cap I_2})\ni S\mapsto 
\big([S]_{I_1/(I_1\cap I_2)},[S]_{I_2/(I_1\cap I_2)}\big)\in {\mathcal M}.$$
is an isomorphism, and its restriction to $\cK(X_{I_1\cap I_2})$ 
is an isomorphism onto the \Csa $\cK$ of ${\mathcal M}$ defined by 
$$\cK=\big\{(k_1,k_2)\in\cK(X_{I_1})\oplus \cK(X_{I_2})\ \big|\ 
[k_1]_{(I_1+I_2)/I_1}=[k_2]_{(I_1+I_2)/I_2}\in \cK(X_{I_1+I_2})\big\}.$$
\end{proposition}

\begin{proof}
Take $(S_1,S_2)\in {\mathcal M}$, 
and we will define $\varPsi'(S_1,S_2)\in \cL(X_{I_1\cap I_2})$. 
For $\xi\in X_{I_1\cap I_2}$, 
we have 
$$[S_1[\xi]_{I_1/(I_1\cap I_2)}]_{(I_1+I_2)/I_1}=
[S_2[\xi]_{I_2/(I_1\cap I_2)}]_{(I_1+I_2)/I_2}.$$
Hence by Proposition \ref{isom3}, 
there exists a unique element $\eta\in X_{I_1\cap I_2}$ 
with 
$$[\eta]_{I_1/(I_1\cap I_2)}=S_1[\xi]_{I_1/(I_1\cap I_2)},\quad
\mbox{and }\quad
[\eta]_{I_2/(I_1\cap I_2)}=S_2[\xi]_{I_2/(I_1\cap I_2)}.$$
We define $\varPsi'(S_1,S_2)\colon X_{I_1\cap I_2}\to X_{I_1\cap I_2}$ 
by $\varPsi'(S_1,S_2)\xi=\eta$ where $\eta$ is the unique element 
satisfying the above two equations. 
It is straightforward to see that 
$$\ip{\varPsi'(S_1,S_2)\xi}{\xi'}_{X_{I_1\cap I_2}}
=\ip{\xi}{\varPsi'(S_1^*,S_2^*)\xi'}_{X_{I_1\cap I_2}}$$
for every $\xi,\xi'\in X_{I_1\cap I_2}$ 
using Lemma \ref{isom1}. 
Thus we have $\varPsi'(S_1,S_2)\in\cL(X_{I_1\cap I_2})$ 
for all $(S_1,S_2)\in {\mathcal M}$. 
It is easy to see that 
$\varPsi'\colon{\mathcal M}\to \cL(X_{I_1\cap I_2})$ 
is a $*$-ho\-mo\-mor\-phism, 
and gives the inverse of $\varPsi$. 
Hence $\varPsi\colon \cL(X_{I_1\cap I_2})\to {\mathcal M}$ 
is an isomorphism. 

Clearly the restriction of $\varPsi$ on $\cK(X_{I_1\cap I_2})$ 
is an injection into $\cK$. 
We will show that this is surjective. 
By Lemma \ref{isom2}, 
we can see that 
the restriction of the map 
$[\cdot]_{I_2/(I_1\cap I_2)}\colon \cK(X_{I_1\cap I_2})\to \cK(X_{I_2})$ to 
$$\ker ([\cdot]_{I_1/(I_1\cap I_2)})=
\cK\big(X_{I_1\cap I_2}\big(I_1/(I_1\cap I_2)\big)\big)$$ 
is a bijection onto 
$$\ker ([\cdot]_{(I_1+I_2)/I_2})=
\cK\big(X_{I_2}\big((I_1+I_2)/I_2\big)\big).$$
Take $(k_1,k_2)\in \cK$. 
Since the map 
$[\cdot]_{I_1/(I_1\cap I_2)}\colon \cK(X_{I_1\cap I_2})\to \cK(X_{I_1})$ 
is surjective, 
we can find $k'\in \cK(X_{I_1\cap I_2})$ 
with $[k']_{I_1/(I_1\cap I_2)}=k_1$. 
Then we see that 
$k_2-[k']_{I_2/(I_1\cap I_2)}\in \ker ([\cdot]_{(I_1+I_2)/I_2})$. 
Thus there exists a unique element 
$k''\in \ker ([\cdot]_{I_1/(I_1\cap I_2)})\subset \cK(X_{I_1\cap I_2})$ 
with $[k'']_{I_2/(I_1\cap I_2)}=k_2-[k']_{I_2/(I_1\cap I_2)}$. 
Now it is easy to see that $k=k'+k''\in \cK(X_{I_1\cap I_2})$ 
satisfies $\varPsi(k)=(k_1,k_2)$. 
We are done. 
\end{proof}

\begin{corollary}\label{isom5}
If $S\in\cL(X_{I_1\cap I_2})$ satisfies 
$$[S]_{I_1/(I_1\cap I_2)}\in\cK(X_{I_1}),\quad 
[S]_{I_2/(I_1\cap I_2)}\in\cK(X_{I_2}),$$
then $S\in\cK(X_{I_1\cap I_2})$. 
\end{corollary}

\begin{proof}
Clear by Proposition \ref{isom4}. 
\end{proof}

\section{\Ccs and representations}\label{Corres}

\begin{definition}\label{DefCor}
For a \Ca $A$, 
we say that $X$ is a {\em $C^*$-cor\-re\-spon\-dence} over $A$ 
when $X$ is a Hilbert $A$-module and 
a $*$-ho\-mo\-mor\-phism $\varphi_X\colon A\to \cL(X)$ is given. 
\end{definition}

We refer to $\varphi_X$ as the left action of 
a \Cc $X$. 
\Ccs can be considered as generalizations of 
automorphisms or endomorphisms. 
In fact, we can associate a \Cc $X_\varphi$ 
with each endomorphism $\varphi$ as follows. 

\begin{definition}
Let $A$ be a $C^*$-al\-ge\-bra 
and $\varphi\colon A\to A$ be an endomorphism. 
We define a \Cc $X_\varphi$ 
such that it is isomorphic to $A$ as Banach spaces, 
its inner product is defined by 
$\ip{\xi}{\eta}_{X}=\xi^*\eta$, 
right action is multiplication 
and left action is given by 
$\varphi_{X_\varphi}(a)\xi=\varphi(a)\xi$. 
We denote $X_{\id_A}$ by $A$, 
and call it the {\em identity correspondence} over $A$. 
\end{definition}

Note that the left action $\varphi_A$ of 
the identity correspondence $A$ 
gives an isomorphism from $A$ to $\cK(A)$. 

\begin{definition}
A {\em morphism} from a \Cc $X$ 
over a \Ca $A$ to 
a \Cc $Y$ over a \Ca $B$ 
is a pair $(\varPi,T)$ consisting of 
a $*$-ho\-mo\-mor\-phism $\varPi\colon A\to B$ and 
a linear map $T\colon X\to Y$ satisfying 
\benu
\item $\bip{T(\xi)}{T(\eta)}_Y=\varPi\big(\ip{\xi}{\eta}_X\big)$ 
for $\xi,\eta\in X$,
\item $\varphi_Y\big(\varPi(a)\big)T(\xi)=
T\big(\varphi_X(a)\xi\big)$ for $a\in A$, $\xi\in X$. 
\eenu
A morphism $(\varPi,T)$ is said to be {\em injective}
if a $*$-ho\-mo\-mor\-phism $\varPi$ is injective. 
\end{definition}

A morphism is called a semicovariant homomorphism in \cite{Sc2}. 
For a morphism $(\varPi,T)$ from $X$ to $Y$, 
we can see that $T(\xi)\varPi(a)=T(\xi a)$ and $\|T(\xi)\|_Y\leq\|\xi\|_X$ 
for $a\in A$ and $\xi\in X$ 
by the same argument as in \cite[Section 2]{Ka5}. 
We also see that $T$ is isometric 
for an injective morphism $(\varPi,T)$. 

\begin{definition}
For a morphism $(\varPi,T)$ from a \Cc $X$ over $A$
to a \Cc $Y$ over $B$, 
we define a $*$-ho\-mo\-mor\-phism $\varPsi_T\colon \cK(X)\to\cK(Y)$ 
by $\varPsi_T(\theta_{\xi,\eta})=\theta_{T(\xi),T(\eta)}$ 
for $\xi,\eta\in X$. 
\end{definition}

For the well-definedness of a $*$-ho\-mo\-mor\-phism $\varPsi_T$, 
see, for example, \cite[Lemma 2.2]{KPW}.
Note that $\varPsi_T$ is injective for an injective morphism $(\varPi,T)$. 
The following two lemmas are easily verified. 

\begin{lemma}\label{psi}
For a morphism $(\varPi,T)$ from a \Cc $X$ over $A$
to a \Cc $Y$ over $B$, 
we have $\varphi_Y(\varPi(a))\varPsi_T(k)=\varPsi_T(\varphi_X(a)k)$ and 
$\varPsi_T(k)T(\xi)=T(k \xi)$ for $a\in A$, $\xi\in X$ and $k\in\cK(X)$. 
\end{lemma}

\begin{lemma}
Let $X$, $Y$, $Z$ be $C^*$-cor\-re\-spon\-dences, 
and $(\varPi_1,T_1)$, $(\varPi_2,T_2)$ be morphisms 
from $X$ to $Y$ and from $Y$ to $Z$, respectively. 
Then its composition $(\varPi_2\circ\varPi_1, T_2\circ T_1)$ is 
a morphism from $X$ to $Z$, 
and 
we have $\varPsi_{T_2\circ T_1}=\varPsi_{T_2}\circ \varPsi_{T_1}$. 
\end{lemma}

\begin{definition}\label{DefRep}
A {\em representation} of a \Cc $X$ over $A$ 
on a \Ca $B$ is 
a pair $(\pi,t)$
consisting of a $*$-ho\-mo\-mor\-phism $\pi\colon A\to B$ 
and a linear map $t\colon X\to B$ satisfying 
\benu
\item $t(\xi)^*t(\eta)=\pi\big(\ip{\xi}{\eta}_X\big)$ 
for $\xi,\eta\in X$,
\item $\pi(a)t(\xi)=t\big(\varphi_X(a)\xi\big)$ for $a\in A$, $\xi\in X$. 
\eenu
We denote by $C^*(\pi,t)$ the \Ca generated 
by the images of $\pi$ and $t$ in $B$. 
We define a $*$-ho\-mo\-mor\-phism $\psi_t\colon \cK(X)\to C^*(\pi,t)$ 
by $\psi_t(\theta_{\xi,\eta})=t(\xi)t(\eta)^*\in C^*(\pi,t)$ 
for $\xi,\eta\in X$. 
\end{definition}

Representations of a \Cc $X$ on a \Ca $B$ is 
nothing but morphisms from $X$ to 
the identity correspondence over $B$, 
and we have $\varphi_B\circ\psi_t=\varPsi_t$. 
Note that we get $\pi(a)\psi_t(k)=\psi_t(\varphi_X(a)k)$ and 
$\psi_t(k)t(\xi)=t(k \xi)$ 
for $k\in\cK(X)$, $a\in A$ and $\xi\in X$. 

\begin{definition}\label{DefGauge}
A representation $(\pi,t)$ of $X$ 
is said to {\em admit a gauge action} 
if for each $z\in\T$, 
there exists a $*$-ho\-mo\-mor\-phism 
$\beta_z\colon C^*(\pi,t)\to C^*(\pi,t)$ 
such that $\beta_z(\pi(a))=\pi(a)$ 
and $\beta_z(t(\xi))=zt(\xi)$ 
for all $a\in A$ and $\xi\in X$. 
\end{definition}

If it exists, such a $*$-ho\-mo\-mor\-phism $\beta_z$ is unique 
and $\beta\colon \T\to\Aut(C^*(\pi,t))$ 
is a strongly continuous homomorphism.

\section{\Cas associated with \Ccs}\label{C*algOfCor}

\begin{definition}\label{DefTX}
For a \Cc $X$ over a \Ca $A$, 
we denote by $\cT_X$ the \Ca generated 
by the universal representation. 
\end{definition}

The universal representation can be obtained by 
taking a direct sum of sufficiently many representations. 
By the universality, 
we have a surjection $\cT_X\to C^*(\pi,t)$ 
for every representation $(\pi,t)$ of $X$. 
The \Ca $\cT_X$ is too large 
to reflect the informations of $X$, 
and so we will take a certain quotient of $\cT_X$ 
to get the nice \Ca $\cO_X$. 

\begin{definition}
For an ideal $I$ of a \Ca $A$, 
we define $I^{\perp}\subset A$ by 
$$I^{\perp}=\{a\in A\mid ab=0 \mbox{ for all }b\in I\}.$$
\end{definition}

Note that 
$I^{\perp}$ is the largest ideal of $A$ 
satisfying $I\cap I^{\perp}=0$. 

\begin{definition}
For a \Cc $X$ over $A$ , 
we define an ideal $J_X$ of $A$ by 
$$J_X=\varphi_X^{-1}\big(\cK(X)\big)\cap \big(\ker\varphi_X\big)^{\perp}.$$
\end{definition}

The ideal $J_X$ is the largest ideal 
to which the restriction of $\varphi_X$ is an injection into $\cK(X)$. 

\begin{definition}\label{DefCov}
A representation $(\pi,t)$ of $X$ is said to be 
{\em covariant} 
if we have $\pi(a)=\psi_t(\varphi_X(a))$ 
for all $a\in J_X$. 
\end{definition}

\begin{definition}\label{DefOX}
For a \Cc $X$ over a \Ca $A$, 
the \Ca $\cO_X$ is defined by $\cO_X=C^*(\pi_X,t_X)$ 
where $(\pi_X,t_X)$ is the universal covariant representation of $X$. 
\end{definition}

By the universality, 
for any covariant representation $(\pi,t)$ of a \Cc $X$, 
there exists a $*$-ho\-mo\-mor\-phism 
$\rho_{(\pi,t)}\colon \cO_X\to C^*(\pi,t)$ 
such that $\pi=\rho_{(\pi,t)}\circ\pi_X$ 
and $t=\rho_{(\pi,t)}\circ t_X$. 
By the universality, 
the universal covariant representation $(\pi_X,t_X)$ 
admits a gauge action. 
We denote it by $\gamma\colon\T\curvearrowright \cO_X$. 
When we consider $\cO_X$ as a generalization 
of crossed products by automorphisms, 
the gauge action $\gamma$ is regarded as 
the dual action of $\T$. 
If a covariant representation $(\pi,t)$ 
admits a gauge action $\beta$, 
then we have $\beta_z\circ\rho_{(\pi,t)}=\rho_{(\pi,t)}\circ\gamma_z$ 
for each $z\in\T$. 
In \cite[Proposition 4.11]{Ka5}, 
we saw that the universal covariant representation $(\pi_X,t_X)$ 
is injective. 
The following gauge-invariant uniqueness theorem says that 
two conditions, admitting a gauge action and being injective, 
characterize the universal one $(\pi_X,t_X)$ 
among all covariant representations. 

\begin{theorem}[{\cite[Theorem 6.4]{Ka5}}]\label{GIUT}
For a covariant representation $(\pi,t)$ of a \Cc $X$, 
the map $\rho_{(\pi,t)}\colon \cO_X\to C^*(\pi,t)$ is an isomorphism 
if and only if $(\pi,t)$ is injective 
and admits a gauge action. 
\end{theorem}

In Proposition \ref{smallest}, 
we see that the universal covariant representation $(\pi_X,t_X)$ 
is the smallest one 
among injective representations admitting gauge actions. 

\begin{remark}
A morphism $(\varPi,T)$ from a \Cc $X$ to a \Cc $Y$ 
gives us a $*$-ho\-mo\-mor\-phism $\cT_X\to\cT_Y$. 
This also gives a $*$-ho\-mo\-mor\-phism $\cO_X\to\cO_Y$ 
when the morphism $(\varPi,T)$ is covariant, that is, 
we have $\varPi(a)\in J_Y$ 
and $\varphi_Y(\varPi(a))=\varPsi_{T}(\varphi_X(a))$
for all $a\in J_X$. 
We do not use these facts explicitly. 
\end{remark}

\section{Invariant ideals}\label{SecInv}

In this section, 
we introduce the notion of invariant ideals 
with respect to $C^*$-cor\-re\-spon\-dences. 
Let us take a \Cc $X$ 
over a \Ca $A$, and fix them until the end of Section \ref{SecIdeal}. 

\begin{definition}
For an ideal $I$ of $A$, 
we define $X(I), X^{-1}(I)\subset A$ by 
\begin{align*}
X(I)&=\cspa\big\{\ip{\eta}{\varphi_X(a)\xi}\in A\ \big|\ 
a\in I,\xi,\eta\in X\big\},\\
X^{-1}(I)&=\{a\in A\mid \ip{\eta}{\varphi_X(a)\xi}_X\in I\mbox{ for all }
\xi,\eta\in X\}.
\end{align*}
\end{definition}

Clearly $X(I)$ is an ideal of $A$. 
We also see that $X^{-1}(I)$ is an ideal 
because it is 
the kernel of the composition of $\varphi_X$ 
and the map $[\cdot]_I\colon \cL(X)\to \cL(X_I)$. 
For a \Cc $X_\varphi$ defined from 
an endomorphism $\varphi\colon A\to A$, 
we see that 
$X_\varphi(I)$ is the ideal generated by $\varphi(I)$, 
and $X_\varphi^{-1}(I)=\varphi^{-1}(I)$ 
for an ideal $I$ of $A$. 
It is easy to see that we have 
$X(I_1)\subset X(I_2)$ and $X^{-1}(I_1)\subset X^{-1}(I_2)$ 
for two ideals $I_1,I_2$ of $A$ with $I_1\subset I_2$. 
For an ideal $I$, 
we have $X(X^{-1}(I))\subset I$ and $X^{-1}(X(I))\supset I$. 
These inclusions are proper in general, 
because we always have $X(I)\subset\cspa\ip{X}{X}_X$ and 
$X^{-1}(I)\supset\ker\varphi_X$. 
The inclusions 
$$X(X^{-1}(I))\subset I\cap \cspa\ip{X}{X}_X,\qquad
X^{-1}(X(I))\supset I+\ker\varphi_X$$
still can be proper as we will see 
in Examples \ref{Ex1} and \ref{Ex2}. 

\begin{lemma}\label{X-1cap}
For two ideals $I_1,I_2$ of $A$, 
we have 
\begin{align*}
X(I_1\cap I_2)&\subset X(I_1)\cap X(I_2),\ \phantom{\mbox{ and }}\ 
X^{-1}(I_1\cap I_2)=X^{-1}(I_1)\cap X^{-1}(I_2),\\
X(I_1+I_2)&=X(I_1)+X(I_2),\ 
\mbox{ and }\ 
X^{-1}(I_1+I_2)\supset X^{-1}(I_1)+X^{-1}(I_2). 
\end{align*}
\end{lemma}

\begin{proof}
Clear by the definitions. 
\end{proof}

Both of the two inclusions in Lemma \ref{X-1cap} can be proper in general
(see Examples \ref{Ex1} and \ref{Ex2}). 

\begin{example}\label{Ex1}
Let $A$ be $\C\oplus\C\oplus M_2(\C)$, 
and $\varphi\colon A\to A$ be an endomorphism defined by 
$\varphi\big((\lambda,\mu,T)\big)=(0,0,\diag\{\lambda,\mu\})$. 
This endomorphism gives us a \Cc $X_\varphi$ over $A$. 
Let us define three ideals $I_1$, $I_2$ and $I_3$ of $A$ by 
$I_1=\C\oplus 0\oplus 0$, $I_2=0\oplus \C\oplus 0$ 
and $I_3=0\oplus 0\oplus M_2(\C)$. 
We see that $\ker\varphi_{X_\varphi}=\ker\varphi=I_3$ 
and $\varphi_{X_\varphi}^{-1}(\cK(X_\varphi))=A$. 
Hence we get $J_X=I_1+I_2$. 
We have $X(I_1)=X(I_2)=I_3$. 
However clearly we have $X(I_1\cap I_2)=X(0)=0$. 
This gives an example of a proper inclusion 
$X(I_1\cap I_2)\subset X(I_1)\cap X(I_2)$. 
Since $X^{-1}(I_3)=A$, 
we have two proper inclusions 
$X^{-1}(X(I_i))\supset I_i+\ker\varphi_{X_\varphi}$ 
for $i=1,2$. 
We see that there exist no non-trivial invariant ideals of $A$ 
(see Definition \ref{Definv}), 
and the \Ca $\cO_X$ is isomorphic to 
a simple \Ca $M_6(\C)$. 
\end{example}

For an increasing family $\{I_n\}_{n\in\N}$ of ideals of a \Ca $D$, 
we denote by $\lim_{n\to\infty}I_n$ the ideal of $D$ defined by 
$$\lim_{n\to\infty}I_n=\overline{\bigcup_{n\in\N}I_n}.$$

\begin{proposition}\label{LimPos}
Let $\{I_n\}_{n\in\N}$ be an increasing family of ideals of $A$. 
Then we have $X(\lim_{n\to\infty}I_n)=\lim_{n\to\infty}X(I_n)$. 
\end{proposition}

\begin{proof}
Clear by the definition of $X(\cdot)$. 
\end{proof}

The analogous statement of Proposition \ref{LimPos} for $X^{-1}$ 
is not valid as the next example shows. 

\begin{example}\label{Ex1.5}
Let $A=C((0,1])$. 
We define a \Cc $X$ over $A$ which is isomorphic to 
$A$ as Hilbert $A$-modules 
and its left action $\varphi_X\colon A\to \cL(X)$ 
is defined by $\varphi_X(f)=f(1)\id_X$ for $f\in A$. 
For each $n\in\N$, 
we define an ideal $I_n$ of $A$ 
by $I_n=C((2^{-n},1])$. 
We have $\lim_{n\to\infty}I_n=A$. 
It is not difficult to see that $X^{-1}(I_n)=C((0,1))$ 
for every $n\in\N$. 
Hence we get $\lim_{n\to\infty}X^{-1}(I_n)=C((0,1))$. 
However, we have $X^{-1}(\lim_{n\to\infty}I_n)=X^{-1}(A)=A$. 
The \Ca $\cO_X$ is isomorphic to the universal \Ca 
generated by a contractive scaling element (see \cite{Ka6}). 
\end{example}

Though we do not have 
$X^{-1}(\lim_{n\to\infty}I_n)=\lim_{n\to\infty}X^{-1}(I_n)$ 
in general, 
we can prove Proposition \ref{LimNeg}, 
which suffices for the further investigation. 
For the proof of Proposition \ref{LimNeg}, 
we need the following general fact. 

\begin{lemma}\label{LemPed}
Let $D$ be a $C^*$-al\-ge\-bra, 
and $\{I_n\}_{n\in\N}$ be 
an increasing family of ideals of $D$. 
For each \Csa $B$ of $D$, 
we have $B\cap (\lim_{n\to\infty}I_n)=\lim_{n\to\infty}(B\cap I_n)$. 
\end{lemma}

\begin{proof}
Set $I_\infty=\lim_{n\to\infty}I_n$. 
Clearly we have 
$B\cap I_\infty \supset\lim_{n\to\infty}(B\cap I_n)$. 
Take a positive element $x\in B\cap I_\infty$. 
For $\e>0$, 
let $f_\e\colon \R_+\to \R_+$ be a continuous function 
defined by $f_\e(t)=\max\{0,t-\e\}$. 
Then we have $\|x-f_\e(x)\|\leq \e$. 
Since $\bigcup_{n\in\N}I_n$ is a dense ideal in $I_\infty$, 
we have $f_\e(x)\in \bigcup_{n\in\N}I_n$ 
(see \cite[Theorem 5.6.1]{Pe}). 
Thus $x$ is approximated by elements 
$f_\e(x)\in B\cap \bigcup_{n\in\N}I_n=\bigcup_{n\in\N}(B\cap I_n)$. 
This shows $x\in \lim_{n\to\infty}(B\cap I_n)$. 
Therefore we have 
$B\cap (\lim_{n\to\infty}I_n)=\lim_{n\to\infty}(B\cap I_n)$. 
\end{proof}

Note that Lemma \ref{LemPed} is not valid 
when $I_n$'s are just $C^*$-sub\-al\-ge\-bras. 

\begin{proposition}\label{LimNeg}
Let $\{I_n\}_{n\in\N}$ be an increasing family of ideals of $A$. 
For each ideal $J$ of $A$ with $\varphi_X(J)\subset \cK(X)$, 
we have 
$J\cap X^{-1}(\lim_{n\to\infty}I_n)=\lim_{n\to\infty}(J\cap X^{-1}(I_n))$. 
\end{proposition}

\begin{proof}
Set $I_\infty=\lim_{n\to\infty}I_n$. 
First note that 
we have 
$$J\cap X^{-1}(I)=\{a\in J\mid \varphi_X(a)\in\cK(XI)\}$$
for an ideal $I$ of $A$ 
by Lemma \ref{ckXI}. 
Take $a\in J\cap X^{-1}(I_\infty)$ and $\e>0$. 
It is easy to see that 
$\cK(XI_\infty)=\lim_{n\to\infty}\cK(XI_n)$. 
By Lemma \ref{LemPed}, 
we have $\varphi_X(J)\cap \cK(XI_\infty) 
=\lim_{n\to\infty}(\varphi_X(J)\cap \cK(XI_n))$. 
Since $\varphi_X(a)\in\varphi_X(J)\cap \cK(XI_\infty)$, 
we can find $n\in\N$ and $k\in \varphi_X(J)\cap \cK(XI_n)$ 
such that $\|\varphi_X(a)-k\|<\e$. 
Then we can find $x\in J$ with $\|x\|<\e$ and $\varphi_X(x)=\varphi_X(a)-k$. 
Set $j=a-x\in J$. 
We have $\varphi_X(j)=k\in \cK(XI_n)$. 
Thus we get $j\in J\cap X^{-1}(I_n)$ and $\|a-j\|<\e$. 
Therefore we get 
$J\cap X^{-1}(I_\infty)
\subset\lim_{n\to\infty}(J\cap X^{-1}(I_n))$. 
The converse inclusion is obvious. 
\end{proof}

\begin{definition}\label{Definv}
An ideal $I$ of $A$ is said to be {\em positively invariant} 
if $X(I)\subset I$, 
{\em negatively invariant} if $J_X\cap X^{-1}(I)\subset I$, 
and {\em invariant} 
if $I$ is both positively and negatively invariant. 
\end{definition}

In many papers such as \cite{KPW}, \cite{FMR} and \cite{Sc2}, 
a positively invariant ideal is called $X$-invariant. 
It is clear that $I$ is positively invariant 
if and only if $I\subset X^{-1}(I)$. 
It is also equivalent to $\varphi_X(I)X\subset XI$ 
by Proposition \ref{XI}. 
Clearly $A$ is an invariant ideal. 
We also see that $0$ is invariant 
because $X(0)=0$ and $J_X\cap X^{-1}(0)=J_X\cap \ker\varphi_X=0$. 

\begin{proposition}\label{LimInv}
Let $\{I_n\}_{n\in\N}$ be an increasing family of ideals of $A$. 
If $I_n$ is positively invariant (resp.\ negatively invariant, invariant), 
then $\lim_{n\to\infty}I_n$ is also. 
\end{proposition}

\begin{proof}
Clear by Proposition \ref{LimPos} and Proposition \ref{LimNeg}. 
\end{proof}

\begin{proposition}\label{capinv}
If two ideals $I_1,I_2$ are 
positively invariant, 
then their intersection $I_1\cap I_2$ is also positively invariant. 
The same is true for negative invariance. 
\end{proposition}

\begin{proof}
Clear by Lemma \ref{X-1cap}. 
\end{proof}

\begin{corollary}\label{capinv2}
The intersection of two 
invariant ideals is invariant. 
\end{corollary}

By Lemma \ref{X-1cap}, 
we see that if two ideals $I_1,I_2$ are 
positively invariant, 
then so is their sum $I_1+I_2$. 
However, the sum of two negatively invariant ideals need not 
be negatively invariant. 
Moreover, the sum of two invariant ideals can fail to be 
negatively invariant as we will see in the next example. 

\begin{example}\label{Ex2}
Let $A$ be $\C\oplus\C\oplus\C$, 
and $X$ be $\C\oplus\C$ which is a Hilbert $A$-module 
by the operations $\ip{(\xi_1,\eta_1)}{(\xi_2,\eta_2)}_{X}
=(\overline{\xi_1}\xi_2,\overline{\eta_1}\eta_2,0)$ 
and $(\xi,\eta)(\lambda,\mu,\nu)
=(\xi\lambda,\eta\mu)$. 
We define a left action $\varphi_X\colon A\to\cL(X)$ 
by $\varphi_X((\lambda,\mu,\nu))=\nu\id_X$. 
We define three ideals $I_1$, $I_2$ and $I_3$ of $A$ by 
$I_1=\C\oplus 0\oplus 0$, $I_2=0\oplus \C\oplus 0$ 
and $I_3=0\oplus 0\oplus \C$. 
We have $J_X=I_3$. 
An easy computation shows that 
$X(I_1)=X(I_2)=0$ and $X^{-1}(I_1)=X^{-1}(I_2)=I_1+I_2$. 
Thus both $I_1$ and $I_2$ are invariant ideals. 
However we have $X(I_1+I_2)=0$ and $X^{-1}(I_1+I_2)=A$. 
Thus $I_1+I_2$ is positively invariant, but not negatively invariant. 
We also have proper inclusions 
$$\begin{array}{rccccclc}
A&=&X^{-1}(I_1+I_2)&\supset& X^{-1}(I_1)+X^{-1}(I_2)&=&I_1+I_2&\\
0&=&X(X^{-1}(I_i))&\subset& I_i\cap \cspa\ip{X}{X}_X&=&I_i&
\quad (i=1,2).
\end{array}$$
We have $\cO_X\cong M_2(\C)\oplus M_2(\C)$, 
and two non-trivial invariant ideals $I_1,I_2$ 
correspond to the two non-trivial ideals of $\cO_X$. 
\end{example}

\begin{definition}
Let us take an ideal $I$ of $A$. 
We define ideals $X^{n}(I)$ for $n\in\N$ 
by $X^0(I)=I$ and $X^{n+1}(I)=X(X^{n}(I))$. 
We also define ideals $X_{-n}(I)$ for $n\in\N$ 
by $X_{0}(I)=I$, 
$X_{-1}(I)=I+J_X\cap X^{-1}(I)$ and 
$X_{-(n+1)}(I)=X_{-1}(X_{-n}(I))$ for $n\geq 1$. 
\end{definition}

Note that we have $I\subset X_{-1}(I)$, 
hence $X_{-n}(I)\subset X_{-(n+1)}(I)$ for every $n\in\N$

\begin{definition}
For an ideal $I$ of $A$, 
we define ideals $X^{\infty}(I)$, $X_{-\infty}(I)$ 
and $X^{\infty}_{-\infty}(I)$ of $A$ by 
$$X^{\infty}(I)=\sum_{n=0}^\infty X^{n}(I)
=\lim_{k\to\infty}\sum_{n=0}^k X^{n}(I),\quad 
X_{-\infty}(I)=\lim_{n\to\infty}X_{-n}(I),$$
$$\mbox{ and }\ X^{\infty}_{-\infty}(I)=X_{-\infty}(X^{\infty}(I)).$$
\end{definition}

\begin{lemma}\label{X-nPos}
If an ideal $I$ is positively invariant, 
so are $X_{-n}(I)$ for $n\in\N\cup\{\infty\}$.
\end{lemma}

\begin{proof}
Let us take a positively invariant ideal $I$. 
From 
$$X_{-1}(I)=I+J_X\cap X^{-1}(I)\subset 
X^{-1}(I)\subset X^{-1}(X_{-1}(I))$$
we see that $X_{-1}(I)$ is positively invariant. 
By using this fact, 
we can prove inductively that $X_{-n}(I)$ is positively invariant 
for all $n\in\N$. 
Finally $X_{-\infty}(I)$ is positively invariant 
by Proposition \ref{LimInv}. 
\end{proof}

\begin{proposition}\label{InvIdealGenBy}
For an ideal $I$ of $A$, 
the ideal $X^{\infty}(I)$ 
(resp.\ $X_{-\infty}(I)$, $X^{\infty}_{-\infty}(I)$) 
is the smallest positively invariant 
(resp.\ negatively invariant, invariant) 
ideal containing $I$. 
\end{proposition}

\begin{proof}
For each $k\in\N$, 
we have 
$$X\bigg(\sum_{n=0}^k X^{n}(I)\bigg)=\sum_{n=0}^k X^{n+1}(I)
\subset X^{\infty}(I).$$
Hence by Proposition \ref{LimPos}, 
we have $X(X^{\infty}(I))\subset X^{\infty}(I)$. 
Thus $X^{\infty}(I)$ is positively invariant. 
If $I'$ is a positively invariant ideal containing $I$, 
then we can prove inductively $X^{n}(I)\subset I'$ for all $n\in\N$. 
Hence we have $X^{\infty}(I)\subset I'$. 
Thus $X^{\infty}(I)$ 
is the smallest positively invariant ideal containing $I$. 

For each $n\in\N$, 
we have 
$J_X\cap X^{-1}(X_{-n}(I))\subset X_{-(n+1)}(I)\subset X_{-\infty}(I)$. 
Hence by Proposition \ref{LimNeg}, 
we have $J_X\cap X^{-1}(X_{-\infty}(I))\subset X_{-\infty}(I)$. 
Thus $X_{-\infty}(I)$ is negatively invariant. 
If $I'$ is a negatively invariant ideal containing $I$, 
then we can prove inductively $X_{-n}(I)\subset I'$ for all $n\in\N$. 
Hence we have $X_{-\infty}(I)\subset I'$. 
Thus $X_{-\infty}(I)$ 
is the smallest negatively invariant ideal containing $I$. 

Combining the above argument 
with Lemma \ref{X-nPos}, 
we see that $X^{\infty}_{-\infty}(I)$ 
is the smallest invariant ideal containing $I$. 
\end{proof}

\section{$T$-pairs and $O$-pairs}\label{SecPairs}

In this section, we introduce the notion of $T$-pairs and $O$-pairs 
of the \Cc $X$ over $A$. 
These are related to 
representations of $X$. 

\begin{definition}
For an ideal $I$ of $A$, 
we define an ideal $J(I)$ of $A$ by 
$$J(I)=\big\{a\in A\ \big|\ [\varphi_X(a)]_I\in\cK(X_{I}),\ 
a X^{-1}(I)\subset I\big\}.$$
\end{definition}

For a positively invariant ideal $I$, 
we can define a map $\varphi_{X_I}\colon A/I\to\cL(X_I)$ 
so that $\varphi_{X_I}([a]_I)=\big[\varphi_{X}(a)\big]_I$ 
because $a\in I$ implies $\big[\varphi_{X}(a)\big]_I=0$. 
Thus in this case, $X_I$ is a \Cc over $A/I$. 
It is clear that the pair $([\cdot]_I,[\cdot]_I)$ of 
the quotient maps $A\to A/I$ and $X\to X_I$ 
is a morphism from $X$ to $X_I$. 

\begin{lemma}\label{JI}
For a positively invariant ideal $I$, 
we have $X^{-1}(I)={[\cdot]_I}^{-1}(\ker\varphi_{X_I})$,
$J(I)={[\cdot]_I}^{-1}(J_{X_I})$ 
and $X^{-1}(I)\cap J(I)=I$. 
\end{lemma}

\begin{proof}
We have 
$$X^{-1}(I)=\ker([\cdot]_I\circ\varphi_X)
=\ker(\varphi_{X_I}\circ[\cdot]_I)={[\cdot]_I}^{-1}(\ker\varphi_{X_I}).$$
We also see that $[\varphi_X(a)]_I\in\cK(X_{I})$ 
if and only if $\varphi_{X_I}([a]_I)\in\cK(X_{I})$. 
Since $X^{-1}(I)={[\cdot]_I}^{-1}(\ker\varphi_{X_I})$, 
the condition $a X^{-1}(I)\subset I$ for $a\in A$ 
is equivalent to $[a]_I \ker\varphi_{X_I}=0$. 
Hence $a\in J(I)$ if and only if 
$$[a]_I\in \varphi_{X_I}^{-1}\big(\cK(X_{I})\big)\cap 
\big(\ker\varphi_{X_I}\big)^{\perp}=J_{X_I}.$$ 
Thus we get $J(I)={[\cdot]_I}^{-1}(J_{X_I})$. 
Finally we have 
$$X^{-1}(I)\cap J(I)={[\cdot]_I}^{-1}(\ker\varphi_{X_I}\cap J_{X_I})
={[\cdot]_I}^{-1}(0)=I.$$
\end{proof}

Note that Lemma \ref{JI} 
implies that $X^{-1}(I)/I=\ker\varphi_{X_I}$ 
and $J(I)/I=J_{X_I}$ 
for a positively invariant ideal $I$. 
Note also that $X^{-1}(0)=\ker\varphi_X$ and $J(0)=J_X$. 

\begin{proposition}\label{neginv}
An ideal $I$ is negatively invariant 
if and only if $J_X\subset J(I)$. 
\end{proposition}

\begin{proof}
For $a\in J_X$, we have $\varphi_X(a)\in\cK(X)$. 
Hence $[\varphi_X(a)]_I\in\cK(X_{I})$. 
Thus $J_X\subset J(I)$ if and only if $J_XX^{-1}(I)\subset I$. 
This is equivalent to the negative invariance of $I$ 
because $J_XX^{-1}(I)=J_X\cap X^{-1}(I)$. 
\end{proof}

Note that $I_1\subset I_2$ need not imply 
$J(I_1)\subset J(I_2)$ in general as the following example shows. 

\begin{example}[{cf.\ Example \ref{Ex2}}]
Let $A\cong\C^3$ be the \Ca generated 
by three mutually orthogonal projections $p_0,p_1$ and $p_2$. 
Let $X$ be the $\ell^\infty$-direct sum 
of two Hilbert spaces $\C$, whose base is denoted by $s_0$, 
and $\ell^2(\N)$, whose base is denoted by $\{s_k\}_{k=1}^\infty$. 
We define an inner product 
$\ip{\cdot}{\cdot}_X\colon X\times X\to A$ 
by $\ip{s_0}{s_0}_X=p_0$, $\ip{s_k}{s_k}_X=p_1$ for $k=1,2,\ldots$, 
and $\ip{s_k}{s_l}_X=0$ for $k\neq l$. 
The right action of $A$ on $X$ is defined by 
$$s_kp_i=\left\{\begin{array}{ll}
s_0&\mbox{for }k=i=0,\\
s_k&\mbox{for }k\geq 1, i=1,\\
0&\mbox{otherwise}.
\end{array}\right.$$
Then $X$ becomes a Hilbert $A$-module. 
We define a left action $\varphi_X\colon A\to\cL(X)$ by 
$\varphi_X(p_0)=\varphi_X(p_1)=0$, and 
$\varphi_X(p_2)=\id_X$. 
Now we get a \Cc $X$ over $A$. 
This \Cc is defined from the following graph (or its opposite graph); 
$$
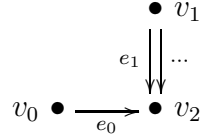
\raisebox{1cm}{\xymatrix{
&\bullet\ v_1 \ar@<-1.8ex>[d]_{e_1} \ar@<-1.0ex>[d]^{\cdots}\\
v_0\ \bullet \ar[r]_{e_0} & \bullet\ v_2}}
$$
(see \cite{Ka2}). 
Let us define ideals of $A$ by 
$$I_0=\C p_0,\ I_1=\C p_1,\ I_{01}=\C p_0+\C p_1\ 
\mbox{ and }\ I_{12}=\C p_1+\C p_2.$$ 
Since $\ker\varphi_X=\varphi_X^{-1}(\cK(X))=I_{01}$, 
we have $J_X=0$. 
Hence all ideals are negatively invariant. 
Since $X(I_1)=X(I_{01})=0$, 
both $I_1$ and $I_{01}$ are invariant. 
By straightforward computation, 
we get $J(I_1)=I_{12}$ and $J(I_{01})=I_{01}$. 
Thus two ideals $I_1,I_{01}$ satisfy that $I_1\subset I_{01}$ 
and $J(I_1)\not\subset J(I_{01})$. 
We can see that $\cO_X$ is isomorphic to 
the direct sum of $M_2(\C)$ and the unitization $\widetilde{K}$ 
of the \Ca $K$ of compact operators on $\ell^2(\N)$. 
There exist six $O$-pairs (see Definition \ref{DefOpair}) 
which correspond to 
six ideals of $\cO_X\cong M_2(\C)\oplus \widetilde{K}$; 
$$
\begin{array}{ccccc}
(0,0)&\subset& (I_1,I_1)&\subset& (I_1,I_{12}) \\
\cap& &\cap& &\cap \\
(I_0,I_0)&\subset& (I_{01},I_{01})&\subset&(A,A)
\end{array}
\ \leftrightsquigarrow\ 
\begin{array}{ccccc}
0 &\subset& K &\subset& \widetilde{K}\\
\cap& &\cap& &\cap \\
M_2(\C) &\subset& M_2(\C)\oplus K &\subset& \cO_X.  
\end{array}
$$
\end{example}

This example also shows that 
$J(I_1\cap I_2)\subset J(I_1)\cap J(I_2)$ 
does not hold in general for two ideals $I_1,I_2$ of $A$. 
However, the converse inclusion 
$J(I_1)\cap J(I_2)\subset J(I_1\cap I_2)$ 
always holds. 

\begin{proposition}\label{capJI}
For two ideals $I_1,I_2$ of $A$, 
we have $J(I_1)\cap J(I_2)\subset J(I_1\cap I_2)$. 
\end{proposition}

\begin{proof}
Take $a\in J(I_1)\cap J(I_2)$. 
Since $[\varphi_X(a)]_{I_1}\in\cK(X_{I_1})$ 
and $[\varphi_X(a)]_{I_2}\in\cK(X_{I_2})$, 
we have $[\varphi_X(a)]_{I_1\cap I_2}\in\cK(X_{I_1\cap I_2})$ 
by Corollary \ref{isom5}. 
We get $a X^{-1}(I_1\cap I_2)\subset I_1\cap I_2$ 
from 
$$a X^{-1}(I_1\cap I_2)\subset a X^{-1}(I_1)\subset I_1,\quad
a X^{-1}(I_1\cap I_2)\subset a X^{-1}(I_2)\subset I_2.$$
Hence $a\in J(I_1\cap I_2)$. 
Thus we have $J(I_1)\cap J(I_2)\subset J(I_1\cap I_2)$. 
\end{proof}

\begin{definition}
Let $X$ be a \Cc over a \Ca $A$. 
A {\em $T$-pair} of $X$ 
is a pair $\omega=(I,I')$ of ideals $I,I'$ of $A$ 
such that $I$ is positively invariant 
and $I\subset I'\subset J(I)$. 
\end{definition}

\begin{definition}
Let $\omega_1=(I_1,I_1')$ and $\omega_2=(I_2,I_2')$
be $T$-pairs. 
We write $\omega_1\subset \omega_2$ 
if $I_1\subset I_2$ and $I_1'\subset I_2'$. 
We denote by $\omega_1\cap \omega_2$ 
the pair $(I_1\cap I_2, I_1'\cap I_2')$. 
\end{definition}

\begin{proposition}
For two $T$-pairs 
$\omega_1=(I_1,I_1'), \omega_2=(I_2,I_2')$, 
their intersection $\omega_1\cap \omega_2=(I_1\cap I_2, I_1'\cap I_2')$ 
is a $T$-pair. 
\end{proposition}

\begin{proof}
By Proposition \ref{capinv}, 
$I_1\cap I_2$ is a positively invariant ideal. 
By Proposition \ref{capJI}, 
we have 
$$I_1\cap I_2\subset I_1'\cap I_2'\subset 
J(I_1)\cap J(I_2)\subset J(I_1\cap I_2).$$
Hence $\omega_1\cap \omega_2$ is a $T$-pair. 
\end{proof}

$T$-pairs arise from representations. 

\begin{definition}\label{DefOmegaRep}
For a representation $(\pi,t)$ of $X$, 
we define $I_{(\pi,t)},I_{(\pi,t)}'\subset A$ by 
$$I_{(\pi,t)}=\ker\pi,\quad 
I_{(\pi,t)}'=\pi^{-1}\big(\psi_t(\cK(X))\big).$$
The pair $\big(I_{(\pi,t)},I_{(\pi,t)}'\big)$ is denoted 
by $\omega_{(\pi,t)}$. 
\end{definition}

Clearly $I_{(\pi,t)}$ is an ideal of $A$. 
By the remark before Definition \ref{DefGauge}, 
we see that $I_{(\pi,t)}'$ is also an ideal of $A$. 

\begin{lemma}\label{II'}
For a representation $(\pi,t)$ 
of a \Cc $X$ over a \Ca $A$, 
we have the following.
\benu
\item $I_{(\pi,t)}$ is positively invariant. 
\item $\ker t=XI_{(\pi,t)}$. 
\item There exists an injective representation $(\dot{\pi},\dot{t})$ 
of the \Cc $X_{I_{(\pi,t)}}$ 
on $C^*(\pi,t)$ such that 
$(\pi,t)=(\dot{\pi}\circ [\cdot]_{I_{(\pi,t)}},
\dot{t}\circ [\cdot]_{I_{(\pi,t)}})$.
\item $a\in I_{(\pi,t)}'$ 
implies $[\varphi_X(a)]_{I_{(\pi,t)}}\in\cK(X_{I_{(\pi,t)}})$ 
and $\pi(a)=\psi_{\dot{t}}([\varphi_X(a)]_{I_{(\pi,t)}})$. 
\item For an element $a\in A$ with $\varphi_X(a)\in\cK(X)$, 
we have $\pi(a)=\psi_t(\varphi_X(a))$ 
if and only if $a\in I_{(\pi,t)}'$. 
\eenu
\end{lemma}

\begin{proof}
\benu
\item 
For $a\in I_{(\pi,t)}$ and $\xi,\eta\in X$, 
we have $\ip{\eta}{\varphi_X(a)\xi}_X\in I_{(\pi,t)}$ because 
$$\pi\big(\ip{\eta}{\varphi_X(a)\xi}_X\big)
=t(\eta)^*t(\varphi_X(a)\xi)
=t(\eta)^*\pi(a)t(\xi)=0.$$
Hence $X(I_{(\pi,t)})\subset I_{(\pi,t)}$. 
Thus $I_{(\pi,t)}$ is positively invariant. 
\item For $\xi\in X$, we have 
\begin{align*}
\xi\in \ker t&\iff t(\xi)=0\iff t(\xi)^*t(\xi)=0\\
&\iff \pi(\ip{\xi}{\xi}_X)=0\iff \ip{\xi}{\xi}_X\in I_{(\pi,t)}
 \iff \xi\in XI_{(\pi,t)}.
\end{align*}
\item Obvious by the definition of $I_{(\pi,t)}$ and (ii). 
\item Since $a\in I_{(\pi,t)}'$, 
we can find $k\in\cK(X)$ with $\pi(a)=\psi_t(k)$. 
For $\xi\in X$, 
we have 
$$t(\varphi_X(a)\xi)=\pi(a)t(\xi)=\psi_t(k)t(\xi)=t(k\xi).$$
Hence 
$\big(\varphi_X(a)-k\big)\xi\in\ker t=XI_{(\pi,t)}$ for all $\xi\in X$. 
This implies that $[\varphi_X(a)]_{I_{(\pi,t)}}=[k]_{I_{(\pi,t)}}\in\cK(X_{I_{(\pi,t)}})$ and 
$$\pi(a)=\psi_t(k)=\psi_{\dot{t}}([k]_{I_{(\pi,t)}})
=\psi_{\dot{t}}([\varphi_X(a)]_{I_{(\pi,t)}}).$$
\item 
If $\pi(a)=\psi_t(\varphi_X(a))$, 
then $a\in I_{(\pi,t)}'$. 
For $a\in I_{(\pi,t)}'$ with $\varphi_X(a)\in\cK(X)$, 
we have 
$\pi(a)
=\psi_{\dot{t}}([\varphi_X(a)]_{I_{(\pi,t)}})=\psi_t(\varphi_X(a))$ 
by (iv). 
\eenu
\end{proof}

\begin{proposition}\label{admpairofrep}
For a representation $(\pi,t)$ of $X$, 
the pair $\omega_{(\pi,t)}$ is a $T$-pair. 
\end{proposition}

\begin{proof}
By Lemma \ref{II'} (i), 
$I_{(\pi,t)}$ is positively invariant. 
Clearly we have $I_{(\pi,t)}\subset I_{(\pi,t)}'$. 
Take $a\in I_{(\pi,t)}'$. 
We have $[\varphi_X(a)]_{I_{(\pi,t)}}\in\cK(X_{I_{(\pi,t)}})$ 
by Lemma \ref{II'} (iv). 
Take $b\in X^{-1}(I_{(\pi,t)})$. 
Since $ab\in I_{(\pi,t)}'$, 
we have $\pi(ab)
=\psi_{\dot{t}}\big([\varphi_X(ab)]_{I_{(\pi,t)}}\big)$
by Lemma \ref{II'} (iv). 
We see $[\varphi_X(ab)]_{I_{(\pi,t)}}=0$ 
because $ab \in X^{-1}(I_{(\pi,t)})$. 
Hence $\pi(ab)=0$. 
Thus we get $ab\in\ker \pi=I_{(\pi,t)}$. 
This shows $a\in J(I_{(\pi,t)})$. 
Hence we get $I_{(\pi,t)}'\subset J(I_{(\pi,t)})$. 
Thus $\omega_{(\pi,t)}=\big(I_{(\pi,t)},I_{(\pi,t)}'\big)$ is 
a $T$-pair. 
\end{proof}

We will see that every $T$-pairs 
come from representations (Proposition \ref{Pairs}). 
By the same way of the proof of Proposition \ref{admpairofrep}, 
we can see that 
for a morphism $(\varPi,T)$ from a \Cc $X$ 
to a \Cc $Y$, 
the pair $\omega_{(\varPi,T)}=\big(I_{(\varPi,T)},I_{(\varPi,T)}'\big)$ 
defined by 
$$I_{(\varPi,T)}=\ker \varPi,\quad 
I_{(\varPi,T)}'=(\varphi_Y\circ\varPi)^{-1}\big(\varPsi_T(\cK(X))\big)$$
is a $T$-pair. 

\begin{definition}\label{DefOpair}
A $T$-pair $\omega=(I,I')$ 
satisfying $J_X\subset I'$ 
is called an {\em $O$-pair}. 
\end{definition}

It is clear that the intersection $\omega_1\cap \omega_2$ 
of two $O$-pairs $\omega_1, \omega_2$ is an $O$-pair. 

\begin{lemma}
A pair $\omega=(I,I')$ of ideals of $A$ 
is an $O$-pair if and only if 
$I$ is invariant and 
$I+J_X\subset I'\subset J(I)$. 
\end{lemma}

\begin{proof}
For an $O$-pair $\omega=(I,I')$, 
we have $I+J_X\subset I'\subset J(I)$. 
Thus we get $J_X\subset J(I)$. 
Now Proposition \ref{neginv} implies that $I$ is negatively invariant. 
Therefore $I$ is an invariant ideal. 
The converse is obvious. 
\end{proof}

For a \Cc $X=C_d(E^1)$ arising from 
a topological graph $E$, 
an $O$-pair $(I,I')$ 
is in the form $(C_0(E^0\setminus X^0),C_0(E^0\setminus Z))$ 
where $(X^0,Z)$ is an admissible pair of closed sets of $E^0$ 
defined in \cite{Ka3}. 

\begin{proposition}\label{cov=O}
A representation $(\pi,t)$ is covariant 
if and only if the pair 
$\omega_{(\pi,t)}$ is an $O$-pair. 
\end{proposition}

\begin{proof}
If $(\pi,t)$ is covariant, 
then clearly $J_X\subset I_{(\pi,t)}'$. 
Thus $\omega_{(\pi,t)}$ is an $O$-pair. 
Conversely, if $\omega_{(\pi,t)}$ is an $O$-pair, 
then for $a\in J_X\subset I_{(\pi,t)}'$, 
we have $\pi(a)=\psi_t(\varphi_X(a))$ 
by Lemma \ref{II'} (v). 
Hence $(\pi,t)$ is covariant. 
\end{proof}

By Proposition \ref{cov=O}, 
we have $\omega_{(\pi,t)}=(0,J_X)$ 
for all injective covariant representations $(\pi,t)$.

\section{\Ccs associated with $T$-pairs}\label{SecC*Ad}

Take a $T$-pair $\omega=(I,I')$ of $X$ 
and fix it throughout this section. 
In this section, 
we construct a \Ca $A_\omega$, 
a \Cc $X_\omega$ over $A_\omega$ 
and a representation $(\pi_\omega,t_\omega)$ 
of $X$ on the \Ca $\cO_{X_\omega}$. 
In the next section, we will see that 
this representation $(\pi_\omega,t_\omega)$ 
has a universal property. 
\begin{center}
\mbox{\xymatrix{
(A,X)\ar_{([\cdot]_I,[\cdot]_I)}[dr]
\ar@{..>}^{(\pi_\omega,t_\omega)}[rrr]&&&
\cO_{X_\omega}\\
&(A/I,X_I)\ar@{=}[d]
\ar^{(\varPi_\omega,T_\omega)}[r]
&(A_\omega,X_\omega)
\ar^{(\varPi,T)}[dl]
\ar_{(\pi_{X_\omega},t_{X_\omega})}[ur]\\
&(A/I,X_I)
}}
\end{center}

\begin{definition}
For a $T$-pair $\omega=(I,I')$ 
of a \Cc $X$ over $A$, 
we define a \Ca $A_\omega$ 
and a Hilbert $A_\omega$-module $X_\omega$ by 
\begin{align*}
A_\omega&=\big\{(b,b')\in A/I\oplus A/I'\ \big|\ 
[b]_{J(I)/I}=[b']_{J(I)/I'}\in A/J(I)\big\},\\
X_\omega&=\big\{(\eta,\eta')\in X_I\oplus X_{I'}\ \big|\ 
[\eta]_{J(I)/I}=[\eta']_{J(I)/I'}\in X_{J(I)}\big\},
\end{align*}
where the operations are defined 
as in Section \ref{HilbMod}. 

\end{definition}

Note that $A_\omega$ is a pull-back \Ca 
of two surjections $[\cdot]_{J(I)/I}\colon A/I\to A/J(I)$ 
and $[\cdot]_{J(I)/I'}\colon A/I'\to A/J(I)$. 

\begin{definition}
We define a $*$-ho\-mo\-mor\-phism 
$\varPsi_\omega\colon \cL(X_I)\to\cL(X_\omega)$ by 
$$\varPsi_\omega(S)(\eta,\eta')=(S\eta,[S]_{I'/I}\eta')\in X_\omega$$
for $S\in \cL(X_I)$ and $(\eta,\eta')\in X_\omega$. 
\end{definition}

\begin{definition}
We define a left action $\varphi_{X_\omega}\colon A_\omega\to\cL(X_\omega)$ by 
$$\varphi_{X_\omega}\big((b,b')\big)
=\varPsi_\omega\big(\varphi_{X_I}(b)\big),$$ 
for $(b,b')\in A_\omega$. 
Thus $X_\omega$ is a \Cc over $A_\omega$. 
\end{definition}

\begin{definition}
We set 
$$\varPi_\omega\colon A/I\ni b\mapsto (b,[b]_{I'/I})\in A_\omega,\quad
T_\omega\colon X_I\ni \eta\mapsto (\eta,[\eta]_{I'/I})\in X_\omega.$$
\end{definition}

\begin{lemma}\label{psisig}
We have 
$\varphi_{X_\omega}\circ\varPi_\omega=\varPsi_\omega\circ\varphi_{X_I}$, 
and $T_\omega(S\eta)=\varPsi_\omega(S)T_\omega(\eta)$
for $S\in \cL(X_I)$ and $\eta\in X_I$. 
\end{lemma}

\begin{proof}
Clear by the definitions.
\end{proof}

From this lemma, 
we easily get the following. 

\begin{proposition}
The pair $(\varPi_\omega,T_\omega)$ 
is an injective morphism from $X_I$ to $X_\omega$, 
and the map $\varPsi_{T_\omega}\colon \cK(X_I)\to\cK(X_\omega)$ 
coincides with the restriction of $\varPsi_\omega$ to $\cK(X_I)$. 
\end{proposition}

The next proposition is also easy to see from the definitions. 

\begin{proposition}
For a $T$-pair $\omega=(I,I')$ with $I'=J(I)$, 
the morphism $(\varPi_\omega,T_\omega)$ 
from $X_I$ to $X_\omega$ is an isomorphism. 
\end{proposition}

To compute $J_{X_\omega}\subset A_\omega$, 
we need the following lemma. 

\begin{lemma}\label{psiS}
A pair $(\varPi,T)$ of maps defined by 
$$\varPi\colon A_\omega\ni (b,b')\mapsto b\in A/I,\quad 
T\colon X_\omega\ni (\eta,\eta')\mapsto\eta\in X_I,$$
is a morphism from $X_\omega$ to $X_I$ satisfying 
$\varPi\circ \varPi_\omega=\id_{A/I}$ and 
$T\circ T_\omega=\id_{X_I}$. 
A $*$-ho\-mo\-mor\-phism 
$\varPsi\colon \cL(X_\omega)\ni S\mapsto T\circ S\circ T_\omega\in \cL(X_I)$ 
satisfies that $\varPsi\circ \varPsi_\omega=\id_{\cL(X_I)}$ 
and the restriction of $\varPsi$ to $\cK(X_\omega)$ 
coincides with $\varPsi_{T}\colon \cK(X_\omega)\to\cK(X_I)$. 
\end{lemma}

\begin{proof}
It is clear that $(\varPi,T)$ is a morphism 
satisfying $\varPi\circ \varPi_\omega=\id_{A/I}$ and 
$T\circ T_\omega=\id_{X_I}$. 
By Lemma \ref{psisig}, 
we have 
$$\varPsi(\varPsi_\omega(S))\eta
=T\big(\varPsi_\omega(S)T_\omega(\eta)\big)
=T\big(T_\omega(S\eta)\big)
=S\eta,$$
for $S\in \cL(X_I)$ and $\eta\in X_I$. 
This proves $\varPsi\circ \varPsi_\omega=\id_{\cL(X_I)}$. 
For $(\eta_1,\eta_1'),(\eta_2,\eta_2')\in X_\omega$ and $\eta\in X_I$, 
we have 
\begin{align*}
\varPsi\big(\theta_{(\eta_1,\eta_1'),(\eta_2,\eta_2')}\big)\eta
&=T\big(\theta_{(\eta_1,\eta_1'),(\eta_2,\eta_2')}T_\omega(\eta)\big)\\
&=T\big((\eta_1\ip{\eta_2}{\eta}_{X_I},
\eta_1'\ip{\eta_2'}{[\eta]_{I'/I}}_{X_{I'}})\big)\\
&=\eta_1\ip{\eta_2}{\eta}_{X_I}\\
&=\theta_{\eta_1,\eta_2}(\eta).
\end{align*}
Hence we have 
$\varPsi\big(\theta_{(\eta_1,\eta_1'),(\eta_2,\eta_2')}\big)
=\theta_{\eta_1,\eta_2}$. 
This shows that the restriction of $\varPsi$ to $\cK(X_\omega)$ 
coincides with $\varPsi_{T}\colon \cK(X_\omega)\to\cK(X_I)$. 
\end{proof}

\begin{proposition}\label{JXomega}
We have 
\begin{align*}
\ker\varphi_{X_\omega}
&=\big\{(b,b')\in A_\omega\ \big|\ b\in\ker\varphi_{X_I}\big\},\\
\varphi_{X_\omega}^{-1}\big(\cK(X_\omega)\big)
&=\big\{(b,b')\in A_\omega\ \big|\ 
b\in\varphi_{X_I}^{-1}\big(\cK(X_I)\big)\big\},\\
J_{X_\omega}
&=\big\{(b,b')\in A_\omega\ \big|\ b\in J_{X_I},\ b'=0\big\}. 
\end{align*}
\end{proposition}

\begin{proof}
Since $\varPsi\circ \varPsi_\omega=\id_{\cL(X_I)}$ by Lemma \ref{psiS}, 
we have $\varPsi\big(\varphi_{X_\omega}((b,b'))\big)=\varphi_{X_I}(b)$ 
for $(b,b')\in A_\omega$. 
Hence for $(b,b')\in A_\omega$, 
we have that $\varphi_{X_\omega}((b,b'))=0$ 
if and only if $\varphi_{X_I}(b)=0$. 
This proves the first equality. 
The second one follows similarly 
because we have 
$\varPsi_\omega\big(\cK(X_I)\big)\subset \cK(X_\omega)$ 
and $\varPsi\big(\cK(X_\omega)\big)\subset\cK(X_I)$. 
We will prove the third equality. 
It is easy to see that for $b\in J_{X_I}$, 
we have 
$$(b,0)\in \varphi_{X_\omega}^{-1}\big(\cK(X_\omega)\big)\cap 
(\ker\varphi_{X_\omega})^{\perp}=J_{X_\omega}.$$ 
Take $(b,b')\in J_{X_\omega}$, 
and we will prove that $b\in J_{X_I}$ and $b'=0$. 
Since $\varphi_{X_\omega}((b,b'))\in\cK(X_\omega)$, 
we have $\varphi_{X_I}(b)\in\cK(X_I)$. 
For any $b_0\in \ker\varphi_{X_I}\subset A/I$, 
we have $\varPi_\omega(b_0)=(b_0,[b_0]_{I'/I})\in \ker\varphi_{X_\omega}$. 
Hence $(b,b')(b_0,[b_0]_{I'/I})=0$. 
This implies that $b\in (\ker\varphi_{X_I})^{\perp}$. 
Hence $b\in J_{X_I}$. 
Since $J_{X_I}=J(I)/I$ by Lemma \ref{JI}, 
we have $[b']_{J(I)/I'}=[b]_{J(I)/I}=0$. 
Hence $(0,{b'}^*)\in A_\omega$. 
Since $(0,{b'}^*)\in \ker\varphi_{X_\omega}$, 
we have $(b,b')(0,{b'}^*)=0$. 
This implies $b'=0$. 
Thus we get 
$J_{X_\omega}=\{(b,b')\in A_\omega\mid b\in J_{X_I},\ b'=0\}$. 
\end{proof}

Note that we have $J_{X_\omega}=\big\{(b,b')\in A_\omega\ \big|\  b'=0\}$
because for $b\in A/I$ we have 
$(b,0)\in A_\omega$ if and only if $b\in J_{X_I}$. 

\begin{definition}\label{DefRepOmega}
We define a $*$-ho\-mo\-mor\-phism $\pi_\omega\colon A\to \cO_{X_\omega}$ 
and a linear map $t_\omega\colon X\to \cO_{X_\omega}$ by 
$$\pi_\omega(a)=\pi_{X_\omega}\big(\varPi_\omega([a]_I)\big),
\quad t_\omega(\xi)=t_{X_\omega}\big(T_\omega([\xi]_I)\big)$$ 
for $a\in A$ and $\xi\in X$, 
where $(\pi_{X_\omega},t_{X_\omega})$ 
is the universal covariant representation 
of the \Cc $X_\omega$ on $\cO_{X_\omega}$. 
\end{definition}

\begin{proposition}\label{surj}
The pair $(\pi_\omega,t_\omega)$ 
is a representation of $X$ on $\cO_{X_\omega}$, 
which admits a gauge action 
and satisfies $C^*(\pi_\omega,t_\omega)=\cO_{X_\omega}$. 
\end{proposition}

\begin{proof}
Since $(\pi_\omega,t_\omega)$ is 
a composition of morphisms, 
it is a representation. 
Clearly the gauge action of $\cO_{X_\omega}$ 
gives a gauge action for 
the representation $(\pi_\omega,t_\omega)$. 
We will prove $C^*(\pi_\omega,t_\omega)=\cO_{X_\omega}$. 
Since $\cO_{X_\omega}$ is 
generated by the images of $\pi_{X_\omega}$ and $t_{X_\omega}$, 
it suffices to show that 
$$\pi_{X_\omega}(A_\omega), t_{X_\omega}(X_\omega)
\subset C^*(\pi_\omega,t_\omega).$$
Take $(b,b')\in A_\omega$. 
Choose $a\in A$ with $[a]_{I'}=b'$. 
We have $b-[a]_{I}\in J(I)/I=J_{X_I}$ 
because 
$[b]_{J(I)/I}=[b']_{J(I)/I'}=[a]_{J(I)}$. 
Thus we have $\varphi_{X_I}(b-[a]_I)\in\cK(X_I)$. 
Hence there exists $k\in\cK(X)$ such that 
$[k]_I=\varphi_{X_I}(b-[a]_I)$. 
Since $(b-[a]_{I},0)\in J_{X_\omega}$ 
by Proposition \ref{JXomega}, 
we have 
\begin{align*}
\pi_{X_\omega}\big((b-[a]_I,0)\big)
&=\psi_{t_{X_\omega}}\big(\varphi_{X_\omega}\big((b-[a]_I,0)\big)\big)\\
&=\psi_{t_{X_\omega}}\big(\varPsi_{\omega}(\varphi_{X_I}(b-[a]_I))\big)\\
&=\psi_{t_{X_\omega}}\big(\varPsi_{t_{\omega}}([k]_I)\big)\\
&=\psi_{t_\omega}(k). 
\end{align*}
Therefore we get 
\begin{align*}
\pi_{X_\omega}\big((b,b')\big)
&=\pi_{X_\omega}\big(([a]_{I},[a]_{I'})\big)
  +\pi_{X_\omega}\big((b-[a]_I,0)\big)\\
&=\pi_\omega(a)+\psi_{t_\omega}(k)\in C^*(\pi_\omega,t_\omega). 
\end{align*}
Thus we have shown that 
$\pi_{X_\omega}(A_\omega)\subset C^*(\pi_\omega,t_\omega)$. 

Take $(\eta,\eta')\in X_\omega$. 
Choose $\xi\in X$ with $[\xi]_{I'}=\eta'$. 
Similarly as above, we get $\eta-[\xi]_I\in X_IJ_{X_I}$. 
Choose $\xi'\in X$ and $b\in J_{X_I}$ with $\eta-[\xi]_I=[\xi']_Ib$. 
Then we have $(\eta-[\xi]_I,0)=T_\omega([\xi']_I) (b,0)$. 
Hence we get 
\begin{align*}
t_{X_\omega}\big((\eta,\eta')\big)
&=t_{X_\omega}\big(([\xi]_I,[\xi]_{I'})\big)
  +t_{X_\omega}\big((\eta-[\xi]_I,0)\big)\\
&=t_\omega(\xi)+
t_{X_\omega}\big(T_\omega([\xi']_I)\big)\pi_{X_\omega}\big((b,0)\big)\\
&=t_\omega(\xi)+
t_\omega(\xi')\pi_{X_\omega}\big((b,0)\big)
\in C^*(\pi_\omega,t_\omega), 
\end{align*}
because $\pi_{X_\omega}\big((b,0)\big)
\in C^*(\pi_\omega,t_\omega)$ as shown above. 
This completes the proof.
\end{proof}

\begin{proposition}\label{Pairs}
For a $T$-pair $\omega=(I,I')$, 
we have $\omega_{(\pi_\omega,t_\omega)}=\omega$. 
\end{proposition}

\begin{proof}
Since the maps $\varPi_\omega\colon A_I\to A_\omega$ 
and $\pi_{X_\omega}\colon A_\omega\to\cO_{X_\omega}$ are injective, 
we have 
$$I_{(\pi_\omega,t_\omega)}=\ker \pi_\omega=\ker ([\cdot]_I)=I.$$ 
For $a\in I'$, 
we have $[a]_I\in I'/I\subset J(I)/I=J_{X_I}$. 
Since $\varPi_\omega([a]_I)=([a]_I,0)\in J_{X_I}$, 
we have 
$$\pi_\omega(a)=\pi_{X_\omega}\big(\varPi_\omega([a]_I)\big)
=\psi_{t_{X_\omega}}\big(\varphi_{X_\omega}(([a]_I,0))\big).$$
We see $\varphi_{X_I}([a]_I)\in\cK(X_I)$ from $[a]_I\in J_{X_I}$. 
Hence by the definition of $\varphi_{X_\omega}$ we get 
$$\varphi_{X_\omega}\big(([a]_I,0)\big)
=\varPsi_{T_\omega}\big(\varphi_{X_I}([a]_I)\big)
\in \varPsi_{T_\omega}\big(\cK(X_I)\big).$$
Since $\cK(X_I)=[\cK(X)]_I$, 
we have 
$$\pi_\omega(a)\in
\psi_{t_{X_\omega}}\big(\varPsi_{T_\omega}\big([\cK(X)]_I\big)\big)
=\psi_{t_\omega}\big(\cK(X)\big).$$
Hence $a\in I_{(\pi_\omega,t_\omega)}'$. 
We have shown that $I'\subset I_{(\pi_\omega,t_\omega)}'$. 
Conversely take $a\in I_{(\pi_\omega,t_\omega)}'$. 
Since 
$$\pi_{X_\omega}\big(\varPi_\omega([a]_I)\big)
=\pi_\omega(a)\in\psi_{t_\omega}\big(\cK(X)\big)\subset 
\psi_{t_{X_\omega}}\big(\cK(X_\omega)\big),$$
we have $\varPi_\omega([a]_I)\in J_{X_\omega}$. 
Hence by Proposition \ref{JXomega}, we have $[a]_{I'}=0$. 
This means $a\in I'$. 
Thus we get $I_{(\pi_\omega,t_\omega)}'\subset I'$. 
Therefore $I_{(\pi_\omega,t_\omega)}'=I'$. 
We have shown that $\omega_{(\pi_\omega,t_\omega)}=\omega$. 
\end{proof}

By Proposition \ref{Pairs}, 
we see that every $T$-pairs come from representations.

\section{\Cas generated by representations of \Ccs}\label{SecC^*Rep}

In this section, we prove the following theorem. 

\begin{theorem}\label{thmfactor}
Let $X$ be a \Cc over a \Ca $A$, 
and $(\pi,t)$ be a representation of $X$. 
If a $T$-pair $\omega$ of $X$ 
satisfies $\omega\subset\omega_{(\pi,t)}$,
then there exists a unique surjective $*$-ho\-mo\-mor\-phism 
$\rho\colon \cO_{X_{\omega}}\to C^*(\pi,t)$ 
such that $\pi=\rho\circ\pi_\omega$ 
and $t=\rho\circ t_\omega$. 
The surjection $\rho$ is an isomorphism 
if and only if $\omega=\omega_{(\pi,t)}$ 
and $(\pi,t)$ admits a gauge action. 
\end{theorem}

Take a representation $(\pi,t)$ 
of a \Cc $X$ 
and a $T$-pair $\omega=(I,I')$ of $X$ 
satisfying $\omega\subset\omega_{(\pi,t)}$. 
In order to get a $*$-ho\-mo\-mor\-phism 
$\rho\colon \cO_{X_{\omega}}\to C^*(\pi,t)$, 
we will construct 
a covariant representation $(\ti{\pi},\tilde{t})$ 
of the \Cc $X_\omega$ on $C^*(\pi,t)$. 
Since $I\subset I_{(\pi,t)}=\ker\pi$, 
we can define a representation $(\dot{\pi},\dot{t})$ 
of a \Cc $X_I$ over $A/I$ on $C^*(\pi,t)$ 
such that $\dot{\pi}([a]_I)=\pi(a)$ for $a\in A$
and $\dot{t}([\xi]_I)=t(\xi)$ for $\xi\in X$ 
as in Lemma \ref{II'} (iii). 
It is easy to see that $I_{(\dot{\pi},\dot{t})}=I_{(\pi,t)}/I$ 
and $I_{(\dot{\pi},\dot{t})}'=I_{(\pi,t)}'/I$. 

\begin{center}
\mbox{\xymatrix{
(A,X)
\ar^{([\cdot]_I,[\cdot]_I)}[dr]
\ar^{(\pi_\omega,t_\omega)}[rrr]
\ar@/_2.5pc/_{(\pi,t)}[ddrr]&&&
\cO_{X_\omega}
\ar@{.>}@/^1pc/^{\rho}[ddl]\\
&(A/I,X_I)\ar_{(\dot{\pi},\dot{t})}[dr]
\ar^{(\varPi_\omega,T_\omega)}[r]
&(A_\omega,X_\omega)
\ar^{(\ti{\pi},\tilde{t})}[d]
\ar^{(\pi_{X_\omega},t_{X_\omega})}[ur]\\
&&C^*(\pi,t)
}}
\end{center}

\begin{definition}
Let $(b,b')\in A_\omega$. 
Take $d\in A/I$ with $[d]_{I'/I}=b'$. 
Define $\ti{\pi}\big((b,b')\big)\in C^*(\pi,t)$ by 
$$\ti{\pi}\big((b,b')\big)
=\dot{\pi}(d)
+\psi_{\dot{t}}\big(\varphi_{X_I}(b-d)\big)\in C^*(\pi,t).$$ 
\end{definition}

Note that this definition makes sense 
because $b-d\in J(I)/I=J_{X_I}$ implies $\varphi_{X_I}(b-d)\in\cK(X_I)$. 
Note also that $\ti{\pi}\big((b,b')\big)\in C^*(\pi,t)$ 
does not depend on the choice of 
$d\in A/I$ with $[d]_{I'/I}=b'$ because we have 
$\dot{\pi}(d_1-d_2)=\psi_{\dot{t}}\big(\varphi_{X_I}(d_1-d_2)\big)$ 
if $d_1-d_2\in I'/I\subset I_{(\pi,t)}'/I=I_{(\dot{\pi},\dot{t})}'$ 
by Lemma \ref{II'} (v).

\begin{lemma}
The map $\ti{\pi}\colon A_\omega\to C^*(\pi,t)$ 
is a $*$-ho\-mo\-mor\-phism. 
\end{lemma}

\begin{proof}
It is obvious that $\ti{\pi}$ 
is a $*$-preserving linear map. 
We will show $\ti{\pi}$ is multiplicative. 
Take $(b_1,b_1'),(b_2,b_2')\in A_\omega$. 
Take $d_1,d_2\in A/I$ 
with $[d_1]_{I'/I}=b_1'$, $[d_2]_{I'/I}=b_2'$. 
Since 
$$\dot{\pi}(d)\psi_{\dot{t}}(\varphi_{X_I}(b))
=\psi_{\dot{t}}(\varphi_{X_I}(db))$$ 
for $d\in A/I$ and $b\in J(I)/I=J_{X_I}$, 
we have 
\begin{align*}
\ti{\pi}\big((&b_1,b_1')\big)\ti{\pi}\big((b_2,b_2')\big)\\
&=\big(\dot{\pi}(d_1)+\psi_{\dot{t}}(\varphi_{X_I}(b_1-d_1))\big)
  \big(\dot{\pi}(d_2)+\psi_{\dot{t}}(\varphi_{X_I}(b_2-d_2))\big)\\
&=\dot{\pi}(d_1d_2)
  +\psi_{\dot{t}}\big(\varphi_{X_I}\big(d_1(b_2-d_2)
  +(b_1-d_1)d_2+(b_1-d_1)(b_2-d_2)\big)\big)\\
&=\dot{\pi}(d_1d_2)
  +\psi_{\dot{t}}\big(\varphi_{X_I}(b_1b_2-d_1d_2)\big)\\
&=\ti{\pi}\big((b_1b_2,b_1'b_2')\big)\\
&=\ti{\pi}\big((b_1,b_1')(b_2,b_2')\big).
\end{align*}
Hence $\ti{\pi}$ is a $*$-ho\-mo\-mor\-phism. 
\end{proof}

\begin{proposition}\label{inj}
The map $\ti{\pi}\colon  A_\omega\to C^*(\pi,t)$ 
is injective 
if and only if $\omega=\omega_{(\pi,t)}$. 
\end{proposition}

\begin{proof}
Suppose that $\ti{\pi}$ is injective. 
For $a\in I_{(\pi,t)}$, 
we have $([a]_I,[a]_{I'})\in A_\omega$ and 
$$\ti{\pi}\big(([a]_I,[a]_{I'})\big)
=\dot{\pi}([a]_I)=\pi(a)=0.$$ 
Hence $([a]_I,[a']_I)=0$. 
This implies $a\in I$. 
Thus we get $I_{(\pi,t)}=I$. 
For $a\in I_{(\pi,t)}'$, 
we have 
$[a]_I\in I_{(\pi,t)}'/I\subset J(I_{(\pi,t)})/I=J(I)/I$. 
Hence we get $(0,[a]_{I'})\in A_\omega$. 
We also get $\varphi_{X_I}([a]_I)\in\cK(X_I)$. 
Since $[a]_I\in I_{(\pi,t)}'/I=I_{(\dot{\pi},\dot{t})}'$, 
we have 
$$\ti{\pi}\big((0,[a]_{I'})\big)
=\dot{\pi}([a]_I)-\psi_{\dot{t}}(\varphi_{X_I}([a]_I))=0,$$ 
by Lemma \ref{II'} (v). 
Since $\ti{\pi}$ is injective, 
we have $(0,[a]_{I'})=0$. 
This implies $a \in I'$. 
Thus we get $I_{(\pi,t)}'=I'$. 
Therefore if $\ti{\pi}$ is injective, 
then $\omega=\omega_{(\pi,t)}$. 

Conversely assume $\omega=\omega_{(\pi,t)}$. 
Take $(b,b')\in A_\omega$ with 
$\ti{\pi}\big((b,b')\big)=0$. 
Take $d\in A/I$ with $[d]_{I'/I}=b'$. 
Then we have $\dot{\pi}(d)=\psi_{\dot{t}}(\varphi_{X_I}(d-b))$. 
Hence $d\in I_{(\pi,t)}'/I=I'/I$. 
Therefore we have $b'=0$. 
We also have $\psi_{\dot{t}}(\varphi_{X_I}(b))=0$. 
Since $I=I_{(\pi,t)}$, 
the map $\dot{t}$ is injective. 
Hence $\psi_{\dot{t}}$ is also injective. 
Therefore we have $b\in\ker \varphi_{X_I}$. 
We also have $b\in J(I)/I=J_{X_I}$ 
because $[b]_{J(I)/I}=[b']_{J(I)/I'}=0$. 
Hence $b=0$. 
We have proved that $\ti{\pi}$ is injective. 
\end{proof}

\begin{definition}
Let $\zeta\in X_IJ_{X_I}$. 
Take $\eta\in X_I$ and $b\in J_{X_I}$ 
such that $\zeta=\eta b$. 
We define 
$\bar{t}(\zeta)=\dot{t}(\eta)\psi_{\dot{t}}(\varphi_{X_I}(b))
\in C^*(\pi,t)$. 
\end{definition}

\begin{lemma}\label{lem4}
The map $\bar{t}\colon X_IJ_{X_I}\to C^*(\pi,t)$ 
is a well-defined linear map 
satisfying that 
$\dot{t}(\eta)^*\bar{t}(\zeta)
=\psi_{\dot{t}}\big(\varphi_X(\ip{\eta}{\zeta}_{X_I})\big)$ 
for all $\zeta\in X_IJ_{X_I}$ and $\eta\in X_I$, and 
$\bar{t}(\zeta_1)^*\bar{t}(\zeta_2)
=\psi_{\dot{t}}\big(\varphi_{X_I}(\ip{\zeta_1}{\zeta_2}_{X_I})\big)$ 
for all $\zeta_1,\zeta_2\in X_IJ_{X_I}$. 
\end{lemma}

\begin{proof}
Take $\eta_1,\eta_2\in X_I$, $b_1,b_2\in J_{X_I}$, 
and define $\zeta_1,\zeta_2\in X_IJ_{X_I}$ by 
$\zeta_1=\eta_1 b_1,\zeta_2=\eta_2 b_2$. 
We have 
\begin{align*}
\dot{t}(\eta_1)^*\dot{t}(\eta_2)\psi_{\dot{t}}(\varphi_{X_I}(b_2))
&=\dot{\pi}(\ip{\eta_1}{\eta_2}_{X_I})\psi_{\dot{t}}(\varphi_{X_I}(b_2))\\
&=\psi_{\dot{t}}\big(\varphi_{X_I}(\ip{\eta_1}{\eta_2}_{X_I}b_2)\big)\\
&=\psi_{\dot{t}}\big(\varphi_{X_I}(\ip{\eta_1}{\zeta_2}_{X_I})\big). 
\end{align*}
Similar computation shows that 
$$\big(\dot{t}(\eta_1)\psi_{\dot{t}}(\varphi_{X_I}(b_1))\big)^*
\big(\dot{t}(\eta_2)\psi_{\dot{t}}(\varphi_{X_I}(b_2))\big)
=\psi_{\dot{t}}\big(\varphi_{X_I}(\ip{\zeta_1}{\zeta_2}_{X_I})\big).$$ 
For $\zeta\in X_IJ_{X_I}$, 
take $\eta_1,\eta_2\in X_I$ and $b_1,b_2\in J_{X_I}$ 
such that $\zeta=\eta_1b_1=\eta_2b_2$. 
Set $x=\dot{t}(\eta_1)\psi_{\dot{t}}(\varphi_{X_I}(b_1))
-\dot{t}(\eta_2)\psi_{\dot{t}}(\varphi_{X_I}(b_2))\in C^*(\pi,t)$. 
We have $x^*x=0$ because for $i,j=1,2$, 
$$\big(\dot{t}(\zeta_i)\psi_{\dot{t}}(\varphi_{X_I}(b_i))\big)^*
\big(\dot{t}(\zeta_j)\psi_{\dot{t}}(\varphi_{X_I}(b_j))\big)
=\psi_{\dot{t}}\big(\varphi_{X_I}(\ip{\zeta}{\zeta}_{X_I})\big).$$ 
This shows $\bar{t}$ is well-defined. 
We can check the linearity of $\bar{t}$ 
in a similar fashion. 
The two equalities in the statement had been already checked. 
\end{proof}

\begin{lemma}\label{lem3}
We have 
$$\dot{\pi}(b)\bar{t}(\zeta)=\bar{t}(\varphi_{X_I}(b)\zeta),\quad 
\psi_{\dot{t}}(k)\bar{t}(\zeta)=\bar{t}(k\zeta),$$
for $b\in A/I$, $k\in\cK(X_I)$, and $\zeta\in X_IJ_{X_I}$. 
\end{lemma}

\begin{proof}
Take $\eta\in X_I$ and $d\in J_{X_I}$ with $\zeta=\eta d$. 
Then we have 
\begin{align*}
\dot{\pi}(b)\bar{t}(\zeta)
&=\dot{\pi}(b)\dot{t}(\eta)\psi_{\dot{t}}(\varphi_{X_I}(d))
&\psi_{\dot{t}}(k)\bar{t}(\zeta)
&=\psi_{\dot{t}}(k)\dot{t}(\eta)\psi_{\dot{t}}(\varphi_{X_I}(d))\\
&=\dot{t}(\varphi_{X_I}(b)\eta)\psi_{\dot{t}}(\varphi_{X_I}(d))
&&=\dot{t}\big(k\eta\big)\psi_{\dot{t}}(\varphi_{X_I}(d))\\
&=\bar{t}\big((\varphi_{X_I}(b)\eta)d\big)
&&=\bar{t}\big((k\eta)d\big)\\
&=\bar{t}(\varphi_{X_I}(b)\zeta), &&=\bar{t}(k\zeta). 
\end{align*}
\end{proof}

\begin{lemma}\label{lem5}
For $\zeta\in X_I(I'/I)$, 
we have $\bar{t}(\zeta)=\dot{t}(\zeta)$. 
\end{lemma}

\begin{proof}
Choose $\eta\in X_I$ and $b\in I'/I\subset J(I)/I=J_{X_I}$ 
such that $\zeta=\eta b$. 
Since $b\in I'/I\subset I_{(\pi,t)}'/I=I_{(\dot{\pi},\dot{t})}'$, 
we have 
$\dot{\pi}(b)=\psi_{\dot{t}}(\varphi_{X_I}(b))$ 
by Lemma \ref{II'} (v).
Hence, we get 
$$\bar{t}(\zeta)=\dot{t}(\eta)\psi_{\dot{t}}(\varphi_{X_I}(b))
=\dot{t}(\eta)\dot{\pi}(b)=\dot{t}(\eta b)=\dot{t}(\zeta).$$
\end{proof}

\begin{definition}
Let $(\eta,\eta')\in X_\omega$. 
Take $\zeta\in X_I$ such that $[\zeta]_{I'/I}=\eta'$. 
Define $\tilde{t}\big((\eta,\eta')\big)\in C^*(\pi,t)$ by 
$$\tilde{t}\big((\eta,\eta')\big)
=\dot{t}(\zeta)+\bar{t}(\eta-\zeta)\in C^*(\pi,t).$$ 
\end{definition}

Note that $\eta-\zeta\in X_IJ_{X_I}$ and that 
$\tilde{t}\colon X_\omega\to C^*(\pi,t)$ 
is a well-defined linear map by Lemma \ref{lem5}. 

\begin{proposition}\label{factorthrough}
The pair $(\ti{\pi},\tilde{t})$ 
is a representation 
of the \Cc $X_\omega$ on $C^*(\pi,t)$ 
such that $\dot{\pi}=\ti{\pi}\circ \varPi_\omega$ 
and $\dot{t}=\tilde{t}\circ T_\omega$. 
\end{proposition}

\begin{proof}
It is easy to see that 
$\dot{\pi}=\ti{\pi}\circ \varPi_\omega$ 
and $\dot{t}=\tilde{t}\circ T_\omega$.
We will check that 
the pair $(\ti{\pi},\tilde{t})$ satisfies 
the two conditions in Definition \ref{DefRep}. 
Take $(\eta_1,\eta_1'),(\eta_2,\eta_2')\in X_\omega$. 
Choose $\zeta_1,\zeta_2\in X_I$ 
with $[\zeta_1]_{I'/I}=\eta_1'$, $[\zeta_2]_{I'/I}=\eta_2'$. 
By Lemma \ref{lem4}, 
we have 
\begin{align*}
\tilde{t}\big((\eta_1,\eta_1')&\big)^*
\tilde{t}\big((\eta_2,\eta_2')\big)\\
=&\big(\dot{t}(\zeta_1)+\bar{t}(\eta_1-\zeta_1)\big)^*
  \big(\dot{t}(\zeta_2)+\bar{t}(\eta_2-\zeta_2)\big)\\
=&\dot{\pi}\big(\ip{\zeta_1}{\zeta_2}_{X_I}\big)\\
 &+\psi_{\dot{t}}\big(\varphi_X(\ip{\zeta_1}{\eta_2-\zeta_2}_{X_I}
   +\ip{\eta_1-\zeta_1}{\zeta_2}_{X_I}
   +\ip{\eta_1-\zeta_1}{\eta_2-\zeta_2}_{X_I})\big)\\
=&\dot{\pi}\big(\ip{\zeta_1}{\zeta_2}_{X_I}\big)
  +\psi_{\dot{t}}\big(\varphi_X(\ip{\eta_1}{\eta_2}_{X_I}
   -\ip{\zeta_1}{\zeta_2}_{X_I})\big)\\
=&\ti{\pi}\big(\big(\ip{\eta_1}{\eta_2}_{X_I},
      \ip{\eta_1'}{\eta_2'}_{X_I}\big)\big)\\
=&\ti{\pi}\big(\bip{(\eta_1,\eta_1')}{(\eta_2,\eta_2')}_{X_\omega}\big).
\end{align*}
This proves the condition (i) in Definition \ref{DefRep}. 
We will check the condition (ii). 
Take $(b,b')\in A_\omega$ and $(\eta,\eta')\in X_\omega$. 
Choose $d\in A/I$ and $\zeta\in X_I$ 
with $[d]_{I'/I}=b'$ and $[\zeta]_{I'/I}=\eta'$. 
By Lemma \ref{lem3}, 
we have 
\begin{align*}
\ti{\pi}\big((b&,b')\big)
\tilde{t}\big((\eta,\eta')\big)\\
=&\big(\dot{\pi}(d)+\psi_{\dot{t}}(\varphi_{X_I}(b-d))\big)
  \big(\dot{t}(\zeta)+\bar{t}(\eta-\zeta)\big)\\
=&\dot{\pi}(d)\dot{t}(\zeta)
  +\psi_{\dot{t}}(\varphi_{X_I}(b-d))\dot{t}(\zeta)\\
 &+\dot{\pi}(d)\bar{t}(\eta-\zeta)
  +\psi_{\dot{t}}(\varphi_{X_I}(b-d))\bar{t}(\eta-\zeta)\\
=&\dot{t}(\varphi_{X_I}(d)\zeta)+\dot{t}(\varphi_{X_I}(b-d)\zeta)
  +\bar{t}(\varphi_{X_I}(d)(\eta-\zeta))
  +\bar{t}(\varphi_{X_I}(b-d)(\eta-\zeta))\\
=&\dot{t}(\varphi_{X_I}(b)\zeta)+\bar{t}(\varphi_{X_I}(b)(\eta-\zeta)).
\end{align*}
On the other hand, 
we have 
$$\varphi_{X_\omega}((b,b'))(\eta,\eta')
=\big(\varphi_{X_I}(b)\eta,[\varphi_{X_I}(b)]_{I'/I}\eta'\big)
=\big(\varphi_{X_I}(b)\eta,[\varphi_{X_I}(b)\zeta]_{I'/I}\big).$$
Hence we get 
$$\tilde{t}\big(\varphi_{X_\omega}((b,b'))(\eta,\eta')\big)
=\dot{t}\big(\varphi_{X_I}(b)\zeta\big)
+\bar{t}\big(\varphi_{X_I}(b)(\eta-\zeta)\big).$$
Thus we have $\ti{\pi}\big((b,b')\big)
\tilde{t}\big((\eta,\eta')\big)
=\tilde{t}\big(\varphi_{X_\omega}((b,b'))(\eta,\eta')\big)$. 
We are done. 
\end{proof}

\begin{proposition}
The representation $(\ti{\pi},\tilde{t})$ 
is covariant. 
\end{proposition}

\begin{proof}
Take $(b,0)\in J_{X_\omega}$. 
By definition, 
we have $\ti{\pi}((b,0))=\psi_{\dot{t}}(\varphi_{X_I}(b))$. 
Since $\varphi_{X_\omega}((b,0))=\varPsi_{T_\omega}(\varphi_{X_I}(b))$, 
we have 
$$\psi_{\tilde{t}}\big(\varphi_{X_\omega}((b,0))\big)
=\psi_{\tilde{t}}\big(\varPsi_{T_\omega}(\varphi_{X_I}(b))\big)
=\psi_{\tilde{t}\circ T_\omega}(\varphi_{X_I}(b))
=\psi_{\dot{t}}(\varphi_{X_I}(b)).$$
Hence we get 
$\ti{\pi}((b,0))=\psi_{\tilde{t}}\big(\varphi_{X_\omega}((b,0))\big)$ 
for every element $(b,0)\in J_{X_\omega}$. 
This completes the proof. 
\end{proof}

\begin{lemma}\label{LemC*ofRep}
The representation $(\ti{\pi},\tilde{t})$ of $X_{\omega}$ 
is injective if and only if $\omega=\omega_{(\pi,t)}$. 
It admits a gauge action if and only if so does $(\pi,t)$. 
\end{lemma}

\begin{proof}
The first assertion follows from Proposition \ref{inj}. 
If a representation $(\pi,t)$
admits a gauge action $\beta$, 
then $\beta$ is also a gauge action 
for the representation $(\ti{\pi},\tilde{t})$ 
because $\beta_z(\psi_t(k))=\psi_t(k)$ 
for all $k\in\cK(X)$ and $z\in\T$. 
The converse is obvious. 
\end{proof}

Now we are ready to prove the main theorem of 
this section. 

\begin{proof}[Proof of Theorem \ref{thmfactor}]
Define 
$\rho=\rho_{(\ti{\pi},\tilde{t})}\colon 
\cO_{X_{\omega}}\to C^*(\pi,t)$. 
Since $\dot{\pi}=\ti{\pi}\circ \varPi_\omega$ 
and $\dot{t}=\tilde{t}\circ T_\omega$, 
we have $\pi=\rho\circ\pi_\omega$ 
and $t=\rho\circ t_\omega$. 
This implies that $\rho$ is surjective. 
The uniqueness follows from $C^*(\pi_\omega,t_\omega)=\cO_{X_\omega}$ 
which was proved in Proposition \ref{surj}. 
Finally by Lemma \ref{LemC*ofRep} and Theorem \ref{GIUT}, 
$\rho$ is an isomorphism 
if and only if $\omega=\omega_{(\pi,t)}$ 
and $(\pi,t)$ admits a gauge action. 
\end{proof}

\begin{corollary}\label{C*ofRep2}
Let $X$ be a \Cc over a \Ca $A$ 
and $(\pi,t)$ be a representation of $X$ 
which admits a gauge action. 
Then the \Ca $C^*(\pi,t)$ 
is naturally isomorphic to the \Ca 
$\cO_{X_{\omega_{(\pi,t)}}}$. 
\end{corollary}

We finish this section by the next result, 
which gives a characterization of the \Ca $\cO_X$ 
without using $J_X$ and the notion of covariance. 

\begin{proposition}\label{smallest}
If a representation $(\pi,t)$ is injective 
and admits a gauge action, 
then there exists a surjection $\overline{\rho}\colon C^*(\pi,t)\to\cO_X$ 
with $\pi_X=\overline{\rho}\circ\pi$ and $t_X=\overline{\rho}\circ t$. 
\end{proposition}

\begin{proof}
Set $\omega=\omega_{(\pi,t)}=(I_{(\pi,t)},I_{(\pi,t)}')$. 
Since $(\pi,t)$ is injective, 
we have $I_{(\pi,t)}=0$ and $I_{(\pi,t)}'\subset J(0)=J_X$. 
Hence we get $\omega\subset (0,J_X)=\omega_{(\pi_X,t_X)}$. 
Thus by Theorem \ref{thmfactor}, 
there exists a surjective $*$-ho\-mo\-mor\-phism 
$\rho\colon \cO_{X_{\omega}}\to\cO_X$ 
with $\pi_X=\rho\circ\pi_\omega$ and $t_X=\rho\circ t_\omega$. 
Since $(\pi,t)$ admits a gauge action, 
the \Ca $C^*(\pi,t)$ is isomorphic to $\cO_{X_{\omega}}$ 
by Corollary \ref{C*ofRep2}. 
This completes the proof. 
\end{proof}

By Proposition \ref{smallest}, 
we can define $\cO_X$ to be the smallest \Ca 
among \Cas generated by 
injective representations admitting gauge actions. 
Now Theorem \ref{GIUT} tells us that 
the covariance of representations 
characterizes the representation $(\pi_X,t_X)$ 
among injective representations admitting gauge actions.

\section{Structure of gauge-invariant ideals of $\cO_X$}\label{SecIdeal}

We say that an ideal of $\cO_X$ is gauge-invariant 
if it is globally invariant under the gauge action $\gamma$. 
In this section, 
we analyze structure of gauge-invariant ideals of $\cO_X$. 

\begin{definition}
For an ideal $P$ of $\cO_X$, 
we define $I_P,I_P'\subset A$ by 
$$\pi_X(I_P)=\pi_X(A)\cap P,\quad 
\pi_X(I_P')=\pi_X(A)\cap \big(P+\psi_{t_X}(\cK(X))\big).$$
We set $\omega_{P}=(I_P,I_P')$. 
\end{definition}

\begin{proposition}
For an ideal $P$ of $\cO_X$, 
denote by $\sigma_P$ a natural surjection from $\cO_X$ to $\cO_X/P$. 
Then we have 
$\omega_{P}=\omega_{(\sigma_P\circ\pi_X,\sigma_P\circ t_X)}$. 
Hence $\omega_{P}$ is an $O$-pair. 
\end{proposition}

\begin{proof}
Clear by the definitions. 
\end{proof}

\begin{definition}
Let $\omega$ be an $O$-pair of $X$. 
The representation $(\pi_\omega,t_\omega)$ of $X$ 
on $\cO_{X_\omega}$ is covariant 
by Proposition \ref{cov=O} and Proposition \ref{Pairs}. 
Hence there exists a surjection 
$\rho_{(\pi_\omega,t_\omega)}\colon \cO_X\to \cO_{X_\omega}$. 
We define $P_\omega=\ker \rho_{(\pi_\omega,t_\omega)}$. 
\end{definition}

\begin{lemma}
For an $O$-pair $\omega$, 
the ideal $P_\omega$ of $\cO_X$ is gauge-invariant 
and satisfies $\omega_{P_\omega}=\omega$. 
\end{lemma}

\begin{proof}
Clear by the definitions. 
\end{proof}

\begin{proposition}\label{PropGII}
For a gauge-invariant ideal $P$ of $\cO_X$, 
we have $P=P_{\omega_P}$ and $\cO_X/P\cong \cO_{X_{\omega_P}}$. 
\end{proposition}

\begin{proof}
If $P$ is gauge-invariant, 
the representation 
$(\sigma_P\circ\pi_X,\sigma_P\circ t_X)$ 
admits a gauge action, 
where $\sigma_P\colon \cO_X\to\cO_X/P$ is a natural surjection. 
Hence by the definition of $\omega_{P}$ 
and Theorem \ref{thmfactor}, 
we have an isomorphism $\rho\colon \cO_{X_{\omega_P}}\to \cO_X/P$ 
such that $(\rho\circ\pi_{\omega_P},\rho\circ t_{\omega_P})=
(\sigma_P\circ\pi_X,\sigma_P\circ t_X)$. 
Hence $\cO_X/P\cong \cO_{X_{\omega_P}}$ and $P=P_{\omega_P}$. 
\end{proof}

Now we get the following. 

\begin{theorem}\label{GII}
The set of all gauge-invariant ideals of $\cO_X$ 
corresponds bijectively to the set of all $O$-pairs of $X$ 
by $P\mapsto \omega_P$ and $\omega\mapsto P_\omega$. 
These maps preserve inclusions and intersections.
\end{theorem}

In the case that \Ccs are defined from graphs, 
or more generally from topological graphs, 
Theorem \ref{GII} had already been proved 
in \cite{BHRS} or \cite{Ka3}. 

\begin{corollary}[{\cite[Theorem 6.4]{MT}}]\label{InvBij}
If $A=J_X+\ker\varphi_X$, 
then $P\mapsto I_P$ is a bijection 
from the set of all gauge-invariant ideals of $\cO_X$ 
to the set of all invariant ideals of $A$ with respect to $X$. 
\end{corollary}

\begin{proof}
By Theorem \ref{GII} and Lemma \ref{JI}, 
it suffices to show that $J_{X_I}\subset [J_X]_I$ 
for all invariant ideal $I$ of $A$. 
Let $I$ be an invariant ideal. 
Since $A=J_X+\ker\varphi_X$, 
we have $A/I=[J_X]_I+[\ker\varphi_X]_I$. 
Hence we get $([\ker\varphi_X]_I)^{\perp}=[J_X]_I$. 
Since $\ker\varphi_{X_I}\supset [\ker\varphi_X]_I$, 
we obtain 
$$J_{X_I}\subset (\ker\varphi_{X_I})^{\perp}
\subset ([\ker\varphi_X]_I)^{\perp}=[J_X]_I.$$
We are done. 
\end{proof}

Note that the assumption $A=J_X+\ker\varphi_X$ 
is equivalent to the assumption in \cite[Theorem 6.4]{MT}. 
This is also equivalent to say that $A\cong A_1\oplus A_2$ 
and $\varphi_X\colon A\to\cL(X)$ is the composition of 
the natural surjection $A\to A_1$ 
and an embedding $A_1\hookrightarrow\cK(X)$. 
This assumption is not neccesary to have 
that the map $P\mapsto I_P$ is bijective, 
as we will see in Sections \ref{SecGI&Morita} and \ref{SecBimod}. 
We finish this section by the following result 
on the gauge-invariant ideals of $\cT_X$. 

\begin{proposition}\label{Tideal}
The set of all gauge-invariant ideals of $\cT_X$ 
corresponds bijectively to the set of all $T$-pairs of $X$ 
such that inclusions and intersections are preserved. 
\end{proposition}

\begin{proof}
The set of all gauge-invariant ideals of $\cT_X$ 
corresponds bijectively to the ``set'' of all representations of $X$ 
admitting gauge actions if we consider 
two representations $(\pi,t)$ and $(\pi',t')$ are same 
when there exists a (necessarily unique) isomorphism 
$\rho\colon C^*(\pi,t)\to C^*(\pi',t')$ 
such that $\rho\circ\pi=\pi'$ and $\rho\circ t=t'$. 
Under this identification, 
the ``set'' of all representations of $X$ 
admitting gauge actions 
corresponds bijectively to the set of all $T$-pairs of $X$ 
by $(\pi,t)\mapsto\omega_{(\pi,t)}$ 
defined in Definition \ref{DefOmegaRep}, 
and $\omega\mapsto (\pi_\omega,t_\omega)$ 
defined in Definition \ref{DefRepOmega} 
by Proposition \ref{Pairs} and Theorem \ref{thmfactor}. 
This completes the proof. 
\end{proof}

\section{Gauge invariant ideals and strong Morita equivalence.}
\label{SecGI&Morita}

In this section, 
we will prove that each gauge-invariant ideal $P$ of the \Ca $\cO_X$ 
is strongly Morita equivalent to the \Ca $\cO_{Y_P}$ 
for a certain \Cc $Y_P$. 
In the next section, 
we will see that in fact we can find a \Cc $Y'_P$ 
so that $P$ is isomorphic to $\cO_{Y_P'}$. 

For a positively invariant ideal $I$ of $A$, 
we have $\varphi_X(I)X\subset XI$. 
Hence the closed linear subspace 
$Y_I=\varphi_X(I)X$ of $X$ 
is naturally considered as a \Cc over $I$. 

\begin{lemma}\label{LemYI}
For a positively invariant ideal $I$ of $A$, 
we have $\ker\varphi_{Y_I}=I\cap \ker\varphi_X$ and 
$\varphi_{Y_I}^{-1}(\cK(Y_I))=I\cap \varphi_X^{-1}(\cK(X))$. 
\end{lemma}

\begin{proof}
Take $a\in \ker\varphi_{Y_I}$. 
For $\xi\in X$, 
we have $\varphi_X(a)\varphi_X(a^*)\xi=0$ 
because $\varphi_X(a^*)\xi\in Y_I$. 
Hence we have $aa^*\in \ker\varphi_X$. 
Thus we get $a\in I\cap \ker\varphi_X$. 
This shows $\ker\varphi_{Y_I}\subset I\cap \ker\varphi_X$. 
Since the converse inclusion is obvious, 
we get $\ker\varphi_{Y_I}=I\cap \ker\varphi_X$. 

Take $a\in \varphi_{Y_I}^{-1}(\cK(Y_I))$. 
Set $k=\varphi_{Y_I}(a)\in \cK(Y_I)\subset \cK(X)$. 
Since we have $\varphi_X(a)\varphi_X(b)\xi=k\varphi_X(b)\xi$ 
for $b\in I$ and $\xi\in X$, 
we get $\varphi_X(a)\varphi_X(a)^*=k\varphi_X(a)^*$. 
We also get $\varphi_X(a)k^*=kk^*$ because $k\in \cK(Y_I)$. 
Thus $(\varphi_X(a)-k)(\varphi_X(a)-k)^*=0$. 
Hence $\varphi_X(a)=k\in \cK(X)$. 
Since the converse inclusion is obvious, 
we have $\varphi_{Y_I}^{-1}(\cK(Y_I))=I\cap \varphi_X^{-1}(\cK(X))$. 
\end{proof}

\begin{proposition}\label{PropYI}
For a positively invariant ideal $I$ of $A$, 
we have $J_{Y_I}=I\cap J_X$. 
\end{proposition}

\begin{proof}
Since $\ker\varphi_{Y_I}\subset \ker\varphi_X$, 
we have $(\ker\varphi_{Y_I})^{\perp}\supset (\ker\varphi_X)^{\perp}$. 
By Lemma \ref{LemYI}, 
we have $(\ker\varphi_{Y_I})^{\perp}\cap I\cap \ker\varphi_X=0$. 
Hence $(\ker\varphi_{Y_I})^{\perp}\cap I\subset (\ker\varphi_X)^{\perp}$. 
Thus we have 
$(\ker\varphi_{Y_I})^{\perp}\cap I=(\ker\varphi_X)^{\perp}\cap I$. 
From this equality and Lemma \ref{LemYI}, we get 
\begin{align*}
J_{Y_I}
&=\varphi_{Y_I}^{-1}\big(\cK(Y_I)\big)
 \cap \big(\ker\varphi_{Y_I}\big)^{\perp}\\
&=I\cap \varphi_X^{-1}(\cK(X))
 \cap \big(\ker\varphi_{Y_I}\big)^{\perp}\\
&=I\cap \varphi_X^{-1}(\cK(X))
 \cap \big(\ker\varphi_X\big)^{\perp}\\ 
&=I\cap J_X.
\end{align*}
\end{proof}

\begin{proposition}\label{PropHer}
For a positively invariant ideal $I$ of $A$, 
the \Csa generated by $\pi_X(I)$ and $t_X(Y_I)$ 
is isomorphic to $\cO_{Y_I}$, 
and it is the smallest hereditary \Csa in $\cO_X$ 
containing $\pi_X(I)$. 
\end{proposition}

\begin{proof}
Let $B$ be the \Csa of $\cO_X$ generated by $\pi_X(I)$ and $t_X(Y_I)$. 
Clearly the restrictions of $\pi_X$ and $t_X$ to $I$ and $Y_I$ 
give an injective representation $(\pi,t)$ of $Y_I$ on $\cO_X$ 
which admits a gauge action. 
It is also clear that $C^*(\pi,t)=B$. 
By Proposition \ref{PropYI}, 
this representation $(\pi,t)$ is covariant. 
Thus $B$ is isomorphic to $\cO_{Y_I}$ 
by Theorem \ref{GIUT}. 

Since we have $\pi_X(I)t_X(X)\pi_X(I)=t_X(Y_I)\pi_X(I)=t_X(Y_I)$, 
$B$ is contained in the \Csa $\pi_X(I)\cO_X\pi_X(I)$. 
By \cite[Proposition 2.7]{Ka5}, 
$\cO_X$ is the closure of the linear span of elements in the form 
$$t_X(\xi_1)\cdots t_X(\xi_n)\pi_X(a)
t_X(\eta_m)^*\cdots t_X(\eta_1)^*$$
for $a\in A$ and $\xi_k,\eta_l\in X$. 
By using the fact $\pi_X(I)t_X(X)=t_X(Y_I)\pi_X(I)$, 
we can prove by the induction on $n$ 
that $\pi_X(b)t_X(\xi_1)\cdots t_X(\xi_n)\pi_X(a)\in B$ 
for $b\in I$, $a\in A$ and $\xi_k\in X$. 
Hence $\pi_X(I)\cO_X\pi_X(I)$ is contained in $B$. 
Thus we have shown that 
$B=\pi_X(I)\cO_X\pi_X(I)$ 
which is the smallest hereditary \Csa containing $\pi_X(I)$. 
\end{proof}

\begin{proposition}\label{PropIdealGen}
For an ideal $I$ of $A$, 
the ideal of $\cO_X$ generated by $\pi_X(I)$ is 
$P_{\omega}$ where 
$\omega=(X_{-\infty}^\infty(I),X_{-\infty}^\infty(I)+J_X)$. 
\end{proposition}

\begin{proof}
Let $P$ be the ideal of $\cO_X$ generated by $\pi_X(I)$. 
Since $I\subset I_P$ and $I_P$ is invariant, 
we have $X_{-\infty}^\infty(I)\subset I_P$ 
by Proposition \ref{InvIdealGenBy}.
Hence we have $X_{-\infty}^\infty(I)+J_X\subset I_P+J_X \subset I_P'$. 
Therefore we get $\omega\subset (I_P,I_P')=\omega_P$. 
Since $\pi_X(I)\subset \pi_X(X_{-\infty}^\infty(I))\subset P_{\omega}$ 
implies $P\subset P_{\omega}$, 
we have $\omega_P\subset \omega_{P_{\omega}}=\omega$. 
Thus we get $\omega_P=\omega$. 
Since $\pi_X(I)$ is closed under the gauge action, 
the ideal $P$ is gauge-invariant. 
Hence we have $P=P_{\omega_P}=P_\omega$ 
by Proposition \ref{PropGII}. 
\end{proof}

\begin{proposition}\label{PropMorita}
Let $I$ be a positively invariant ideal of $A$. 
For an $O$-pair $\omega=(X_{-\infty}(I),X_{-\infty}(I)+J_X)$, 
the gauge-invariant ideal $P_\omega$ 
is strongly Morita equivalent to the \Ca $\cO_{Y_I}$. 
\end{proposition}

\begin{proof}
By Proposition \ref{PropHer} and Proposition \ref{PropIdealGen}, 
the \Csa generated by $\pi_X(I)$ and $t_X(Y_I)$ 
is isomorphic to $\cO_{Y_I}$ 
and is a hereditary and full \Csa of $P_\omega$ 
which is the ideal generated by $\pi_X(I)$. 
Thus $P_\omega$ is strongly Morita equivalent to the \Ca $\cO_{Y_I}$. 
\end{proof}

\begin{corollary}\label{CorNonDeg}
Let $X$ be a \Cc over a \Ca $A$. 
Define a \Cc $Y$ over $A$ by $Y=\varphi_X(A)X$. 
Then $\cO_Y$ is strongly Morita equivalent to $\cO_X$. 
\end{corollary}

\begin{proof}
Apply Proposition \ref{PropMorita} to the invariant ideal $A$. 
\end{proof}

The \Cc $Y$ defined in the above corollary 
is non-degenerate, namely 
it satisfies that $\varphi_Y(A)Y=Y$. 
Thus by Corollary \ref{CorNonDeg}, 
we can exchange a given \Cc to a non-degenerate one 
so that the \Cas constructed 
by them are strongly Morita equivalent 
(we used this fact in \cite[Appendix C]{Ka5}). 

By Proposition \ref{PropMorita}, 
gauge-invariant ideals $P$ satisfying that $I_P'=I_P+J_X$ 
are shown to be strongly Morita equivalent 
to the \Ca $\cO_{Y_{I_P}}$ of the \Cc $Y_{I_P}$. 
To deal with all gauge-invariant ideals of $\cO(X)$, 
we need the following argument. 

Let us define a \Ca $\widetilde{A}$ 
and a Banach space $\widetilde{X}$ by 
\begin{align*}
\widetilde{A}&=\pi_X(A)+\psi_{t_X}(\cK(X))\subset\cO_X,\\ 
\widetilde{X}&=\cspa \big(t_X(X)+t_X(X)\psi_{t_X}(\cK(X))\big)\subset\cO_X. 
\end{align*}
If we define the left and right actions of $\widetilde{A}$ 
on $\widetilde{X}$ as multiplication, 
and the inner product by $\ip{\xi}{\eta}_{\widetilde{X}}=\xi^*\eta$, 
$\widetilde{X}$ becomes a \Cc over $\widetilde{A}$. 
Since the embeddings $\widetilde{A}\hookrightarrow \cO_X$ and 
$\widetilde{X}\hookrightarrow \cO_X$ 
give an injective representation of $\widetilde{X}$, 
we have an injective $*$-ho\-mo\-mor\-phism from 
$\cK(\widetilde{X})$ onto 
$\cspa (\widetilde{X}\widetilde{X}^*)\subset\cO_X$. 
Thus we can identify $\cK(\widetilde{X})$ 
with $\cspa (\widetilde{X}\widetilde{X}^*)$. 

\begin{lemma}\label{JwX}
We have $J_{\widetilde{X}}=\psi_{t_X}(\cK(X))\subset \widetilde{A}$. 
\end{lemma}

\begin{proof}
By the identification above, 
the restriction of $\varphi_{\widetilde{X}}$ 
to the ideal $\psi_{t_X}(\cK(X))$ of $\widetilde{A}$ 
is just the embedding $\psi_{t_X}(\cK(X))\hookrightarrow \cK(\widetilde{X})$. 
Hence we have $\psi_{t_X}(\cK(X))\subset J_{\widetilde{X}}$. 
We will prove the converse inclusion. 
Take $\pi_X(a)+\psi_{t_X}(k)\in J_{\widetilde{X}}$. 
Then we have $\pi_X(a)\in J_{\widetilde{X}}$. 
Let $\{u_\lambda\}$ be an approximate unit of $\psi_{t_X}(\cK(X))$. 
It is not difficult to see that $\{\varphi_{\widetilde{X}}(u_\lambda)\}$ 
is an approximate unit of $\cK(\widetilde{X})$ 
(see \cite[Lemma 5.10]{Ka5}).
Since $\varphi_{\widetilde{X}}(\pi_X(a))\in \cK(\widetilde{X})$, 
we have 
$$\varphi_{\widetilde{X}}(\pi_X(a))
=\lim_{\lambda}\varphi_{\widetilde{X}}(\pi_X(a))
\varphi_{\widetilde{X}}(u_\lambda)
=\lim_{\lambda}\varphi_{\widetilde{X}}(\pi_X(a)u_\lambda)
\in \varphi_{\widetilde{X}}\big(\psi_{t_X}(\cK(X))\big).$$ 
Hence there exists $k\in \cK(X)$ with 
$\varphi_{\widetilde{X}}(\pi_X(a))=\varphi_{\widetilde{X}}(\psi_{t_X}(k))$. 
Therefore we have 
\begin{align*}
t_X(\varphi_X(a)\xi)=\pi_X(a)t_X(\xi)
&=\varphi_{\widetilde{X}}(\pi_X(a))t_X(\xi)\\
&=\varphi_{\widetilde{X}}(\psi_{t_X}(k))t_X(\xi)
=\psi_{t_X}(k)t_X(\xi)
=t_X(k\xi)
\end{align*}
for each $\xi\in X$. 
Hence we obtain $\varphi_X(a)=k\in \cK(X)$. 
For $b\in\ker\varphi_X$ 
we have $\pi_X(b)\in \ker\varphi_{\widetilde{X}}$. 
Therefore we get $\pi_X(ab)=0$ for all $b\in\ker\varphi_X$. 
Thus $a\in \varphi_X^{-1}(\cK(X))\cap (\ker\varphi_X)^{\perp}=J_X$. 
Therefore 
$\pi_X(a)+\psi_{t_X}(k)=\psi_{t_X}(\varphi_X(a)+k)\in \psi_{t_X}(\cK(X))$. 
This shows $J_{\widetilde{X}}\subset \psi_{t_X}(\cK(X))$. 
Thus we get $J_{\widetilde{X}}=\psi_{t_X}(\cK(X))$. 
\end{proof}

\begin{proposition}
The natural inclusions $\widetilde{A}\hookrightarrow \cO_X$ and 
$\widetilde{X}\hookrightarrow \cO_X$ induce an isomorphism 
$\cO_{\widetilde{X}}\cong \cO_X$. 
\end{proposition}

\begin{proof}
It is clear that the pair $(\pi,t)$ of the inclusions 
$\pi\colon \widetilde{A}\hookrightarrow \cO_X$ and 
$t\colon \widetilde{X}\hookrightarrow \cO_X$ 
is an injective representation of $\widetilde{X}$ 
admitting a gauge action 
and satisfying $C^*(\pi,t)=\cO_X$. 
By Lemma \ref{JwX}, 
the representation $(\pi,t)$ is covariant. 
Hence we have an isomorphism 
$\rho_{(\pi,t)}\colon \cO_{\widetilde{X}}\to \cO_X$ 
by Theorem \ref{GIUT}. 
\end{proof}

\begin{proposition}\label{PropMorita2}
For a gauge-invariant ideal $P$ of $\cO_X$, 
we set $\widetilde{I}=\widetilde{A}\cap P$. 
Then $P$ is strongly Morita equivalent to the \Ca 
$\cO_{Y_{\widetilde{I}}}$ 
where $Y_{\widetilde{I}}
=\varphi_{\widetilde{X}}(\widetilde{I})\widetilde{X}$ 
is a \Cc over $\widetilde{I}$. 
\end{proposition}

\begin{proof}
Since $\widetilde{I}$ is the intersection of $\widetilde{A}$ 
and the ideal $P$ of $\cO_{\widetilde{X}}=\cO_X$, 
the ideal $\widetilde{I}$ is an invariant ideal of $\widetilde{A}$. 
Let $\widetilde{P}$ be the ideal in $\cO_{\widetilde{X}}=\cO_X$ 
generated by $\widetilde{I}$. 
By Proposition \ref{PropMorita}, 
$\widetilde{P}$ is strongly Morita equivalent to the \Ca 
$\cO_{Y_{\widetilde{I}}}$. 
We will show that $\widetilde{P}=P$. 
To do so, 
it suffices to see $\omega_{\widetilde{P}}=\omega_P$ 
by Theorem \ref{GII} 
because both $\widetilde{P}$ and $P$ are gauge-invariant. 
Since $\widetilde{I}\subset P$, 
we have $\widetilde{P}\subset P$. 
Hence $\omega_{\widetilde{P}}\subset \omega_P$. 
We have 
$$\pi_X(A)\cap P=\pi_X(A)\cap \widetilde{A}\cap P
=\pi_X(A)\cap \widetilde{I}\subset \pi_X(A)\cap \widetilde{P}.$$
Similarly, we have 
\begin{align*}
\pi_X(A)\cap \big(P+\psi_{t_X}(\cK(X))\big)
&=\pi_X(A)\cap \big(\widetilde{A}\cap P+\psi_{t_X}(\cK(X))\big)\\
&=\pi_X(A)\cap \big(\widetilde{I}+\psi_{t_X}(\cK(X))\big)\\
&\subset \pi_X(A)\cap \big(\widetilde{P}+\psi_{t_X}(\cK(X))\big).
\end{align*}
Hence we get $\omega_P\subset \omega_{\widetilde{P}}$. 
Thus $\omega_{\widetilde{P}}=\omega_P$. 
This completes the proof. 
\end{proof}

\begin{remark}
As we saw in the proof of Proposition \ref{PropMorita2}, 
we can see that gauge-invariant ideals of $\cO_X$ are 
distinguished by their intersection with $\widetilde{A}$. 
By Proposition \ref{PropMorita2}, 
the set of all gauge-invariant ideals of $\cO_{\widetilde{X}}$ 
corresponds bijectively 
to the set of all invariant ideals of $\widetilde{A}$ 
even though the \Cc $\widetilde{X}$ does not satisfy 
the assumption in Corollary \ref{InvBij} in general. 
\end{remark}

Proposition \ref{PropMorita2} shows that 
every gauge-invariant ideals of $\cO_X$ are 
strongly Morita equivalent to the \Ca $\cO_{Y}$ 
for some \Ccs $Y$. 
In the next section, 
we will see that for every gauge-invariant ideal $P$ of $\cO_X$ 
we can find a \Cc $Y$ so that $P$ is isomorphic to $\cO_{Y}$.

\section{Crossed products by Hilbert $C^*$-bimodules}\label{SecBimod}

For a \Ca $A$, 
a {\em Hilbert $A$-bimodule} is a \Cc $X$ over $A$ together with 
a left inner product $\lip{\cdot}{\cdot}\colon X\times X\to A$ 
such that $\varphi_X(\lip{\xi}{\eta})=\theta_{\xi,\eta}$ 
for $\xi,\eta\in X$ 
(for the detail, see \cite{AEE}, for example). 
We have 
$$J_X=\cspa\{\lip{\xi}{\eta}\in A\mid \xi,\eta\in X\}.$$ 
A \Cc $X$ has a left inner product 
so that it becomes a Hilbert $A$-bimodules 
if and only if we have $\varphi_X(J_X)=\cK(X)$, 
and in this case a left inner product is uniquely 
determined by the structure of \Cc as 
$\lip{\xi}{\eta}=(\varphi_X|_{J_X})^{-1}(\theta_{\xi,\eta})\in J_X$ 
(see \cite[Lemma 3.4]{Ka4}). 

For a general \Cc $X$ over $A$, 
an ideal $I$ of $A$ is positively invariant 
if and only if $\varphi_X(I)X\subset XI$. 
For Hilbert $C^*$-bimodules, 
we get an analogous statement for negative invariance.
Let us fix a Hilbert $A$-bimodule $X$ 
whose left inner product is denoted by $\lip{\cdot}{\cdot}$. 

\begin{lemma}\label{HilbNeg}
An ideal $I$ of $A$ is negatively invariant 
if and only if $\varphi_X(I)X\supset XI$. 
\end{lemma}

\begin{proof}
Let $I$ be a negatively invariant ideal of $A$. 
Take $\xi\in X$ and $a\in I$. 
For arbitrary $\eta\in X$, 
we have $\varphi_X(\lip{\xi a}{\eta})=\theta_{\xi a,\eta}\in\cK(XI)$. 
Since $\lip{\xi a}{\eta}\in J_X$, 
the negative invariance of $I$ implies 
$\lip{\xi a}{\eta}\in I$ for arbitrary $\eta\in X$. 
Similarly to the proof of Proposition \ref{XI}, 
we can prove $\xi a\in \varphi_X(I)X$. 
Thus we have $\varphi_X(I)X\supset XI$. 
Conversely, assume that an ideal $I$ satisfies $\varphi_X(I)X\supset XI$. 
For $\xi,\eta\in XI$, 
we can find $\xi'\in X$ and $a\in I$ with $\xi=\varphi_X(a)\xi'$. 
Therefore we have 
$\lip{\xi}{\eta}=\lip{\varphi_X(a)\xi'}{\eta}=a(\lip{\xi'}{\eta})\in I$. 
Hence we can see that 
$(\varphi_X|_{J_X})^{-1}(k)\in I$ for $k\in\cK(XI)$. 
Therefore for $a\in J_X$ with $\varphi_X(a)\in\cK(XI)$ 
we have $a\in I$. 
This shows that $I$ is negatively invariant. 
\end{proof}

\begin{proposition}\label{HilbInv}
An ideal $I$ of $A$ is invariant 
if and only if $\varphi_X(I)X=XI$. 
\end{proposition}

\begin{proof}
Clear from Lemma \ref{HilbNeg}. 
\end{proof}

\begin{proposition}\label{HilbQu}
For an invariant ideal $I$ of $A$, 
the \Cc $X_I$ defined in Section \ref{SecPairs} 
has a left inner product ${}_{X_I}\langle{\cdot},{\cdot}\rangle$ 
such that 
${}_{X_I}\langle{[\xi]_I},{[\eta]_I}\rangle=[\lip{\xi}{\eta}]_I$ 
for $\xi,\eta\in X$. 
\end{proposition}

\begin{proof}
Since $\varphi_X(I)X=XI$, 
it is not difficult to see 
that the left inner product of $X_I$ 
described above is well-defined, 
and satisfies the conditions required. 
\end{proof}

\begin{corollary}\label{NoI'}
For an invariant ideal $I$ of $A$, 
we have $J_{X_I}=[J_X]_I$. 
\end{corollary}

\begin{proof}
By Proposition \ref{HilbQu}, 
we have 
\begin{align*}
J_{X_I}
&=\cspa\{{}_{X_I}\langle{\xi'},{\eta'}\rangle\in A/I\mid \xi',\eta'\in X_I\}\\
&=\cspa\{[\lip{\xi}{\eta}]_I\in A/I\mid \xi,\eta\in X\}=[J_X]_I.
\end{align*}
\end{proof}

\begin{proposition}\label{HilbIdeal}
For an invariant ideal $I$, 
the \Csa of $\cO_X$ generated 
by $\pi_X(I)$ and $t_X(XI)$ is an ideal. 
\end{proposition}

\begin{proof}
This follows from the fact that $XI=\varphi_X(I)X=\varphi_X(I)XI$. 
\end{proof}

\begin{theorem}\label{HilbGII}
Let $X$ be a Hilbert $A$-bimodule. 
For an ideal $P$ of $\cO_X$, 
we define an ideal $I_P$ of $A$ by $\pi_X(I_P)=\pi_X(A)\cap P$. 
Then the map $P\mapsto I_P$ gives a one-to-one correspondence 
from the set of all gauge-invariant ideals $P$ of $\cO_X$ 
to the set of ideals $I$ of $A$ 
satisfying $\varphi_X(I)X=XI$. 
We also have isomorphisms $P\cong \cO_{XI_P}$ 
and $\cO_X/P\cong \cO_{X_{I_P}}$ 
for a gauge-invariant ideal $P$. 
\end{theorem}

\begin{proof}
By Corollary \ref{NoI'}, 
we have $I'=I+J_X$ for all $O$-pair $\omega=(I,I')$. 
Thus the first assertion follows from Theorem \ref{GII} 
and Proposition \ref{HilbInv}. 
The second assertion follows from 
Proposition \ref{PropGII}, Proposition \ref{PropHer} 
and Proposition \ref{HilbIdeal}. 
\end{proof}

Note that both $XI$ and $X_I$ are Hilbert $C^*$-bimodule. 
Thus the class of \Cas associated with Hilbert $C^*$-bimodule 
behave well. 
We will see that this class is same as 
the one of \Cas associated with $C^*$-cor\-re\-spon\-dences, 
which we are studying in this paper. 

Let us take a $C^*$-al\-ge\-bra $A$ 
and a $C^*$-cor\-re\-spon\-dence $X$ over $A$. 
We define a \Ca $\overline{A}$ and a Banach space $\overline{X}$ by 
$$\overline{A}=\cO_X^{\,\gamma},\quad 
\overline{X}=\{x\in \cO_X\mid \gamma_z(x)=zx\mbox{ for all }z\in\T\}.$$ 

\begin{remark}
In a similar way to the proof of \cite[Proposition 5.7]{Ka5}, 
we can prove that $\overline{X}=\cspa (t_X(X)\cO_X^{\,\gamma})$. 
We do not use this fact. 
\end{remark}

It is easy to see that $\overline{X}$ is 
a Hilbert $\overline{A}$-bimodule 
where the inner products are defined by 
$$\ip{\xi}{\eta}_{\overline{X}}=\xi^*\eta,\quad 
{}_{\overline{X}}\langle{\xi},{\eta}\rangle=\xi\eta^*,$$
for $\xi,\eta\in \overline{X}$, 
and the left and right actions are multiplication. 

\begin{proposition}[{cf. \cite[Theorem 3.1]{AEE}}]\label{Cc=Hilb}
The natural embedding of $\overline{A}$ and $\overline{X}$ 
into $\cO_X$ 
gives an isomorphism $\cO_{\overline{X}}\cong\cO_X$. 
\end{proposition}

\begin{proof}
By Theorem \ref{GIUT}, 
it suffices to check that 
the embedding of $\overline{A}$ and $\overline{X}$ 
into $\cO_X$ is an injective covariant representation 
admitting a gauge action, 
and these are easily checked. 
\end{proof}

\begin{corollary}\label{GII=OX}
Let $X$ be a \Cc over a \Ca $A$, 
and $P$ be a gauge-invariant ideal of $\cO_X$. 
If we set $I=P\cap \overline{A}$, 
then $P$ is isomorphic to $\cO_{\overline{X}I}$. 
\end{corollary}

\begin{proof}
Combine Theorem \ref{HilbGII} and Proposition \ref{Cc=Hilb}. 
\end{proof}

We remark that 
in order to compute the $K$-groups of gauge-invariant ideals, 
Proposition \ref{PropMorita} and 
Proposition \ref{PropMorita2} seem to be more useful 
than Corollary \ref{GII=OX}.

\section{Relative Cuntz-Pimsner algebras}\label{SecRCP}

In the last section, 
we apply the results obtained above 
to the relative Cuntz-Pimsner algebras 
introduced in \cite{MS}. 
Recall that 
for a \Cc $X$ over a \Ca $A$, 
and an ideal $J$ of $A$ with $\varphi_X(J)\subset\cK(X)$, 
the {\em relative Cuntz-Pimsner algebra} $\cO(J,X)$ 
is generated by the image of a representation $(\pi,t)$ 
which is universal among representations satisfying 
$\pi(a)=\psi_t(\varphi_X(a))$ for $a\in J$ 
(see \cite[Theorem 2.19]{MS}). 
We will show that 
every relative Cuntz-Pimsner algebras 
are isomorphic to $\cO_{X'}$ 
for some \Ccs $X'$. 
In particular, 
every Cuntz-Pimsner algebras and Toeplitz algebras 
(and augmented ones) introduced in \cite{Pi} 
are in the class of our $C^*$-al\-ge\-bras. 

By the universality, 
the representation $(\pi,t)$ 
of $X$ on $\cO(J,X)$ 
admits a gauge action. 
Hence by Corollary \ref{C*ofRep2} 
we see that $\cO(J,X)$ is isomorphic to 
$\cO_{X_{\omega_{(\pi,t)}}}$. 
We will express $\omega_{(\pi,t)}$ 
in terms of a \Cc $X$ over $A$ 
and an ideal $J$ of $A$. 

Now let us take
a \Cc $X$ over a \Ca $A$, 
and an ideal $J$ of $A$ with $\varphi_X(J)\subset\cK(X)$. 
We inductively define an increasing family of ideals 
$\{J_{-n}\}_{n\in\N}$ by 
$J_0=0$ and $J_{-(n+1)}=J_{-n}+J\cap X^{-1}(J_{-n})$. 
We set $J_{-\infty}=\lim_{n\to\infty}J_{-n}$. 
We denote by $\omega_{J}$ the pair $(J_{-\infty},J)$ 
of ideals of $A$. 
Since $X^{-1}(0)=\ker\varphi_X$, 
we have $J_{-1}=J\cap \ker\varphi_X$. 
It is easy to see that $J_{-\infty}=0$ 
if and only if $J\cap \ker\varphi_X=0$. 

\begin{lemma}
The pair $\omega_{J}=(J_{-\infty},J)$ 
is a $T$-pair of $X$. 
\end{lemma}

\begin{proof}
Clearly $J_0=0$ is positively invariant. 
We can prove that $J_{-n}$ is positively invariant for all $n\in\N$ 
by the induction with respect to $\N$ 
similarly as in Lemma \ref{X-nPos}. 
Hence $J_{-\infty}$ is positively invariant. 
Again by the induction, we see that $J_{-\infty}\subset J$. 
Since $J\cap X^{-1}(J_{-n})\subset J_{-(n+1)}\subset J_{-\infty}$ for all $n$, 
we have $J\cap X^{-1}(J_{-\infty})\subset J_{-\infty}$ 
by Proposition \ref{LimNeg}. 
Since $\varphi_X(J)\subset \cK(X)$ by the assumption, 
we have $\big[\varphi_X(J)\big]_{J_{-\infty}}\subset \cK(X_{J_{-\infty}})$. 
Hence we get $J\subset J(J_{-\infty})$. 
Thus we have $J_{-\infty}\subset J\subset J(J_{-\infty})$. 
We are done. 
\end{proof}

\begin{lemma}\label{MinAd}
If a $T$-pair $\omega=(I,I')$ 
satisfies $J\subset I'$, 
then $\omega_{J}\subset \omega$. 
\end{lemma}

\begin{proof}
We will prove $J_{-n}\subset I$ by the induction on $n$. 
For $n=0$, it is trivial. 
Assume $J_{-n}\subset I$. 
We have 
$$J\cap X^{-1}(J_{-n})\subset I'\cap X^{-1}(I)
\subset J(I)\cap X^{-1}(I)=I$$ 
by Lemma \ref{JI}. 
Hence 
$$J_{-(n+1)}=J_{-n}+J\cap X^{-1}(J_{-n})\subset I.$$
We have shown that $J_{-n}\subset I$ for all $n$. 
This implies that $J_{-\infty}\subset I$. 
Hence we have $\omega_{J}\subset \omega$. 
\end{proof}

\begin{proposition}\label{RCP=OX}
The relative Cuntz-Pimsner algebra $\cO(J,X)$ 
is isomorphic to the \Ca $\cO_{X_{\omega_{J}}}$ 
of the \Cc $X_{\omega_{J}}$. 
\end{proposition}

\begin{proof}
Let us denote by $(\pi,t)$ the universal representation 
of $X$ on $\cO(J,X)$ satisfying 
$\pi(a)=\psi_t(\varphi_X(a))$ for all $a\in J$, 
and by $(\pi_{\omega_{J}},t_{\omega_{J}})$ 
the representation of $X$ on the \Ca $\cO_{X_{\omega_{J}}}$ 
defined in Section \ref{SecC*Ad}. 
By Proposition \ref{Pairs}, 
we have $I_{(\pi_{\omega_{J}},t_{\omega_{J}})}'=J$. 
Hence by Lemma \ref{II'} (v), 
we have $\pi_{\omega_{J}}(a)=\psi_{t_{\omega_{J}}}(\varphi_X(a))$ 
for all $a\in J$. 
By the universal property of $\cO(J,X)$, 
there exists a $*$-ho\-mo\-mor\-phism $\rho\colon \cO(J,X)\to \cO_{X_{\omega_{J}}}$ 
such that $\pi_{\omega_{J}}=\rho\circ\pi$ 
and $t_{\omega_{J}}=\rho\circ t$. 
On the other hand, 
$J\subset I_{(\pi,t)}'$ implies $\omega_{J}\subset\omega_{(\pi,t)}$ 
by Lemma \ref{MinAd}. 
Hence by Theorem \ref{thmfactor}, 
there exists a surjective $*$-ho\-mo\-mor\-phism 
$\rho'\colon \cO_{X_{\omega_{J}}}\to \cO(J,X)$ 
such that $\pi=\rho'\circ\pi_{\omega_{J}}$ 
and $t=\rho'\circ t_{\omega_{J}}$. 
Clearly $\rho$ and $\rho'$ are the inverses of each others. 
Hence $\cO(J,X)$ 
is isomorphic to $\cO_{X_{\omega_{J}}}$. 
\end{proof}

By Proposition \ref{RCP=OX}, 
the $T$-pair $\omega_{(\pi,t)}$ 
arising from the representation $(\pi,t)$ on $\cO(J,X)$ 
coincides with $\omega_{J}=(J_{-\infty},J)$. 
From this fact, 
we have the following corollaries. 

\begin{corollary}\label{Cor1}
Let $(\pi,t)$ be the representation of $X$ on $\cO(J,X)$. 
Then the kernel of the map $\pi\colon A\to \cO(J,X)$ 
is $J_{-\infty}$, 
and we have 
$$\{a\in A\mid \varphi_X(a)\in\cK(X),\mbox{ and }
\pi(a)=\psi_t(\varphi_X(a))\}=J.$$
\end{corollary}

\begin{proof}
This easily follows from $\omega_{(\pi,t)}=\omega_{J}$. 
\end{proof}

\begin{corollary}
The relative Cuntz-Pimsner algebra $\cO(J,X)$ is zero 
if and only if $J_{-\infty}=A$. 
\end{corollary}

\begin{proof}
Clear by Proposition \ref{RCP=OX}. 
\end{proof}

\begin{corollary}[{\cite[Proposition 2.21]{MS}}]\label{RCPinj}
The map $\pi\colon A\to \cO(J,X)$ 
is injective if and only if $J\cap \ker\varphi_X=0$. 
\end{corollary}

\begin{proof}
By Corollary \ref{Cor1}, 
$\pi\colon A\to \cO(J,X)$ is injective 
if and only if $J_{-\infty}=0$, 
which is equivalent to the condition 
$J\cap \ker\varphi_X=0$ as we saw above. 
\end{proof}

The following is a gauge-invariant uniqueness theorem 
for relative Cuntz-Pimsner algebras. 

\begin{corollary}\label{Cor2}
For a representation $(\pi',t')$ of $X$ 
satisfying $\pi'(a)=\psi_{t'}(\varphi_X(a))$ for $a\in J$, 
the natural surjection $\cO(J,X)\to C^*(\pi',t')$ 
is an isomorphism if and only if $(\pi',t')$ 
admits a gauge action, $\ker \pi'=J_{-\infty}$, 
and 
$$\{a\in A\mid \pi'(a)\in\psi_{t'}(\cK(X))\}=J.$$
\end{corollary}

\begin{proof}
By Proposition \ref{RCP=OX}, 
$\cO(J,X)$ is canonically isomorphic to 
$\cO_{X_{\omega_{J}}}$. 
By Theorem \ref{thmfactor}, 
the surjection from $\cO_{X_{\omega_{J}}}$ to 
$C^*(\pi',t')$ is injective if and only if 
$(\pi',t')$ admits a gauge action and 
$\omega_{(\pi',t')}=\omega_{J}$. 
The last two conditions in the statement 
just rephrase the condition 
$\omega_{(\pi',t')}=\omega_{J}$. 
\end{proof}

Note that we automatically have $\ker \pi'\supset J_{-\infty}$ and 
$\{a\in A\mid \pi'(a)\in\psi_{t'}(\cK(X))\}\supset J$. 
Note also that in general we cannot replace the condition 
$$\{a\in A\mid \pi'(a)\in\psi_{t'}(\cK(X))\}=J.$$
to the condition 
$$\{a\in A\mid \varphi_X(a)\in\cK(X),\mbox{ and }
\pi'(a)=\psi_{t'}(\varphi_X(a))\}=J,$$
which seems to be natural at first glance. 
This is 
because there may exist $a\in A$ with $\varphi_X(a)\notin\cK(X)$ 
satisfying $[\varphi_X(a)]_{J_{-\infty}}\in \cK(X_{J_{-\infty}})$ 
and 
$\pi'(a)=\psi_{\dot{t'}}([\varphi_X(a)]_{J_{-\infty}})\in \psi_{t'}(\cK(X))$ 
(see Lemma \ref{II'} (iv) and (v)). 
In the case that $J\cap \ker\varphi_X=0$, 
the statement of Corollary \ref{Cor2} 
has the following simple forms. 

\begin{corollary}
Let us assume $J\cap \ker\varphi_X=0$. 
For a representation $(\pi',t')$ of $X$ 
satisfying $\pi'(a)=\psi_{t'}(\varphi_X(a))$ for $a\in J$, 
the natural surjection $\cO(J,X)\to C^*(\pi',t')$ 
is an isomorphism if and only if 
$(\pi',t')$ is injective, admits a gauge action, and satisfies 
$$\{a\in A\mid \varphi_X(a)\in\cK(X),\mbox{ and }
\pi'(a)=\psi_{t'}(\varphi_X(a))\}=J.$$
\end{corollary}

We remark that an ideal $J$ of $A$ satisfies 
$\varphi_X(J)\subset\cK(X)$ and $J\cap \ker\varphi_X=0$ 
if and only if $J\subset J_X$. 
As we saw in Corollary \ref{RCPinj}, 
the maps from $A$ and $X$ to 
the relative Cuntz-Pimsner algebra $\cO(J,X)$ 
is injective 
only when $J$ satisfies $J\subset J_X$. 
Thus it is not a good idea to examine the structure of $\cO(J,X)$ 
in terms of $A$, $X$ and $J$ 
unless $J$ satisfies $J\subset J_X$. 
Anyway, the following result 
on the ideal structure of relative Cuntz-Pimsner algebras $\cO(J,X)$ 
can be easily obtained 
similarly as Theorem \ref{GII} or Proposition \ref{Tideal}. 

\begin{proposition}
Let $X$ be a \Cc over a \Ca $A$, 
and $J$ be an ideal of $A$ with $\varphi_X(J)\subset\cK(X)$. 
Then there exists a one-to-one correspondence 
between the set of all gauge-invariant ideals of $\cO(J,X)$ 
and the set of all $T$-pairs 
$\omega=(I,I')$ of $X$ satisfying $J\subset I'$, 
which preserves inclusions and intersections. 
\end{proposition}

We note that a $T$-pair $\omega=(I,I')$ satisfies $J\subset I'$ 
if and only if $\omega_J\subset\omega$ by Lemma \ref{MinAd}.

\end{document}